\documentclass{article}
\usepackage{amssymb,amsmath,amsthm,mathrsfs}

\newcommand{\ie}{\emph{i.e.}}
\newcommand{\eg}{\emph{e.g.}}
\newcommand{\cf}{\emph{cf}}
\newcommand{\etc}{\emph{etc}}
\newcommand{\Com}{\mathbb{C}}
\newcommand{\Real}{\mathbb{R}}
\newcommand{\Nat}{\mathbb{N}}
\newcommand{\NatStar}{\mathbb{N}^*}

\newcommand{\sii}{L^2}
\newcommand{\sinf}{L^\infty}
\newcommand{\Smooth}{C}
\newcommand{\Dom}{D}
\newcommand{\id}{1}
\newcommand{\Hilbert}{\mathcal{H}}
\newcommand{\Bounded}{\mathfrak{B}}
\newcommand{\demi}{\frac{1}{2}}
\newcommand{\Ball}{B}
\newcommand{\dist}{\mathrm{dist}}
\newcommand{\supp}{\mathop{\mathrm{supp}}\nolimits}
\newcommand{\sgn}{\mathop{\mathrm{sgn}}\nolimits}

\newcommand{\diag}{\mathop{\mathrm{diag}}\nolimits}
\newcommand{\eps}{\varepsilon}
\newcommand{\curve}{\Gamma}
\newcommand{\cross}{\omega}
\newcommand{\tubemap}{\mathcal{L}}

\newcommand{\tube}{\mathcal{T}}
\newcommand{\nodal}{\mathcal{N}}

\newcommand{\cyl}{\Omega}
\newcommand{\transef}{\mathcal{J}}
\newcommand{\curv}{\mathcal{K}}
\newcommand{\rot}{\mathcal{R}}
\newcommand{\vol}{\mathrm{vol}}
\newcommand{\ambient}{\mathcal{A}}
\newtheorem{Claim}{Claim}[section]
\newtheorem{Lemma}[Claim]{Lemma}
\newtheorem{Proposition}[Claim]{Proposition}
\newtheorem{Theorem}[Claim]{Theorem}
\newtheorem{Corollary}[Claim]{Corollary}
\newtheorem{Conjecture}{Conjecture}
\theoremstyle{definition}
\newtheorem{Definition}[Claim]{Definition}
\newtheorem{Remark}[Claim]{Remark}

\newtheorem{ass}{Assumption}
\newenvironment{Assumption}{\begin{ass}}{\end{ass}}

\newcommand{\mybibitem}[1]{\vspace{-1.34ex}\bibitem{#1}}
\begin{document}
%
\title{{\Large\textbf{
Location of the nodal set for thin curved tubes
}}}
\author{
\textsc{
P.~Freitas$^1$
\ and \
D.~Krej\v{c}i\v{r}\'{\i}k$^2$
}}
\date{\footnotesize
\begin{quote}
\emph{
\begin{itemize}
\item[$1$]
Department of Mathematics,
Faculdade de Motricidade Humana (TU Lisbon), \emph{and}
Group of Mathematical Physics, University of Lisbon,
Complexo Interdisciplinar,
Av.~Prof.~Gama Pinto~2,
P-1649-003 Lisboa, Portugal
\smallskip \\
\emph{E-mail:} freitas@cii.fc.ul.pt
\smallskip
\item[$2$]
Department of Theoretical Physics, Nuclear Physics Institute, \\
Academy of Sciences,
250\,68 \v{R}e\v{z} near Prague, Czech Republic
\smallskip \\
\emph{E-mail:} krejcirik@ujf.cas.cz
\end{itemize}
}
\end{quote}
11 January 2007
}
\maketitle
%
%
\begin{abstract}
\noindent
The Dirichlet Laplacian in curved tubes
of arbitrary constant cross-section
rotating together with the Tang frame
along a bounded curve in Euclidean spaces of arbitrary dimension
is investigated in the limit when
the volume of the cross-section diminishes.
We show that spectral properties of the Laplacian
are, in this limit, approximated well by those of the sum
of the Dirichlet Laplacian in the cross-section
and a one-dimensional Schr\"odinger operator
whose potential is expressed solely
in terms of the first curvature of the reference curve.
In particular, we establish the convergence of eigenvalues,
the uniform convergence of eigenfunctions
and locate the nodal set of the Dirichlet Laplacian in the tube
near nodal points of the one-dimensional Schr\"odinger operator.
As a consequence, we prove the ``nodal-line conjecture''
for a class of non-convex and possibly multiply connected domains.
The results are based on a perturbation theory
developed for Schr\"odinger-type operators
in a straight tube of diminishing cross-section.
\bigskip
\begin{itemize}
\item[\textbf{MSC\,2000:}]
Primary 35J25; Secondary 35B25, 35B38.
\item[\textbf{Keywords:}]
Dirichlet Laplacian, nodal set, tubes,
convergence of eigenfunctions.
\end{itemize}
\end{abstract}
%
%
\newpage
\setcounter{equation}{0}
\section{Introduction}
%
Consider the Dirichlet eigenvalue problem
for the Laplacian in a bounded domain $U\subset\Real^d$,
$d \geq 2$:
\begin{equation}\label{problem}
  \left\{
  \begin{aligned}
    -\Delta u &= \lambda u
    &\quad\mbox{in}\quad & U \,,
    \\
    u &= 0
    &\quad\mbox{on}\quad & \partial U \,,
  \end{aligned}
  \right.
\end{equation}
and let us arrange the eigenvalues
in a non-decreasing sequence $\{\lambda_n\}_{n=1}^\infty$
with repetitions according to multiplicities.
The set of corresponding eigenfunctions $\{u_n\}_{n=1}^\infty$
may be chosen in such a way that it forms
an orthonormal basis for~$\sii(U)$.
Since the first eigenfunction~$u_1$ does not vanish in~$U$,
all other eigenfunctions must change sign, and it makes
thus sense to introduce the \emph{nodal set} of~$u_n$ ($n\geq 2$):
\begin{equation}
  \nodal(u_n) := u_n^{-1}(0)
  \,.
\end{equation}
The connected components of $U\setminus\nodal(u_n)$
are called \emph{nodal domains} of~$u_n$.

Since the solutions of~(\ref{problem}) are analytic in~$U$,
each $\nodal(u_n)$ decomposes into the disjoint union
of an analytic $(d-1)$-dimensional manifold
and a singular set contained in a countable number
of analytic $(d-2)$-dimensional manifolds
(\cf~\cite{Caffarelli-Friedman_1985}).
The Courant nodal domain theorem then states that,
if the boundary~$\partial U$ is sufficiently regular,
the $n^\mathrm{th}$~eigenfunction~$u_n$
has at most~$n$ nodal domains
(\cf~\cite{CH1,Chavel,Alessandrini_1998}).
In particular, $u_2$ has exactly two nodal domains.

Apart from these basic results, not much is known regarding
the structure of the nodal sets. One direction along which
much work has been developed over the last three decades
centres around a conjecture of Payne's from 1967~\cite{Payne1},
which states that the nodal set of a second eigenfunction
of problem~(\ref{problem}) for any planar domain
cannot consist of a closed curve.
This conjecture can be extended in an obvious way to higher
dimensions as follows:
\begin{Conjecture}\label{Conj.all}
$
  \overline{\nodal(u_2)} \cap \partial U
  \not= \varnothing
$
for any bounded domain $U\subset\Real^d$.
\end{Conjecture}

The most general result obtained so far
was given by Melas in 1992~\cite{Melas},
who showed that Conjecture~\ref{Conj.all} holds in the case
of planar convex domains (\cf~also~\cite{Alessandrini_1994}).
This followed a string of results by several authors under some
additional symmetry restrictions on the convex domain
\cite{Payne2,Lin_1987,Putter,Damascelli}.

Independently of Melas, Jerison had already
announced a proof of Conjecture~\ref{Conj.all} in 1991
\cite{Jerison0}, in the case of planar convex domains
which are sufficiently long and thin.
In spite of the supplementary eccentricity condition,
Jerison's method has some advantages over Melas's in that,
on the one hand, it also applies to higher dimensions --
\cf~\cite{Jerison} --, and, on the other hand,
it is more quantitative, giving some indication
as to the location of the nodal set.
More precisely, in~\cite{Jerison1},
Jerison located the nodal set in the two-dimensional case
near the zero of an ordinary differential equation
associated to the convex domain in a natural way
(\cf~also~\cite{GJerison}).

On the {\it negative} side, several counterexamples
to Conjecture~\ref{Conj.all} have also been presented.
The most significant is that in~\cite{H2ON} showing that the
result does not hold for multiply connected domains in general
(\cf~also~\cite{Fournais}).
Other counterexamples have been given illustrating different ways in which
the conjecture may not hold. These include adding a
potential~\cite{LinNi} and the case of simply-connected
surfaces~\cite{Freitas1}.
The restriction to bounded domains is also crucial,
since there are examples of simply-connected, unbounded, planar
domains which furthermore satisfy the symmetry restrictions under which
the conjecture holds for the bounded case, but for which the nodal
set does not touch the boundary~\cite{FK2}.

It thus remains an open question for which classes of domains
Conjecture~\ref{Conj.all} holds.
One possibility which seems reasonable is the following
\begin{Conjecture}\label{Conj.simply}
Conjecture~\ref{Conj.all} holds for simply-connected domains
in Euclidean space.
\end{Conjecture}

Let us also mention that the study of nodal sets and domains of
eigenfunctions may of course be extended in a natural way to manifolds
\cite{Cheng-76,Berard-Meyer_1982,
Donnelly-Fefferman_1988,Donnelly-Fefferman_1990a,
Donnelly-Fefferman_1990b,Schoen-Yau,Chavel,Freitas1}.
In fact, in the list of open problems
given in~\cite{Schoen-Yau},
Yau asks to what extent this type of
results can be extended to manifolds -- see Problem 45 in the Chapter
{\it Open problems in differential geometry} in~\cite{Schoen-Yau},
for instance.

The main goal of this paper is to support the validity
of Conjecture~\ref{Conj.simply} by showing that
it holds when~$U$ is a sufficiently thin
curved (and therefore non-convex) \emph{tube} in~$\Real^d$,
for any $d \geq 2$. Note that we allow for the tube to
have an arbitrary cross-section
(rotated appropriately
with respect to the Frenet frame
of a reference curve),
and thus we do not exclude the case
of multiply connected domains either.

This result may be extended to higher eigenfunctions and
we show that, given a natural number $N$ greater than or equal to two,
for any $2\leq n\leq N$ there are precisely $n$ nodal domains of~$u_n$,
and the closure of each of these domains has a non-empty intersection
with~$\partial U$,
provided the tube~$U$ is sufficiently thin.
Moreover, we locate the nodal set~$\nodal(u_n)$
near the zeros of the solution of an ordinary differential equation
which is associated to the tube in a natural way,
via the geometry of the reference curve.

Although the nature of our results
and the main idea behind them are somewhat similar
to those of Jerison's paper~\cite{Jerison1},
the technical approach is actually slightly different.
While in that case the location of the nodal set in a convex domain
is based on the usage of a trial function
and refined applications of the maximum principle,
we use the fact that
the eigenvalue problem for the Laplacian in a curved tube
can be transformed into an eigenvalue problem
for a Schr\"odinger-type operator in a straight tube.
This idea goes at least as far back
as the paper of Exner and \v{S}eba's from 1989~\cite{ES},
where it was used to prove the existence
of discrete eigenvalues in infinite tubes
(\cf~also~\cite{DE,KT}).
This will enable us to develop a perturbation theory
(in which the eigenfunctions of a comparison operator
play the role of the trial function of Jerison's)
and prove $\sii$-convergence results.
However, the maximum principle (in addition to other techniques)
will be also useful for us eventually,
in order to obtain the necessary $\Smooth^0$-convergence results.
More importantly, the main difference
with respect to Jerison's paper
lies in the different setting:
while in Jerison's paper the cross-section
varies wildly along a straight line
(and part of the point is to show that
under the convexity assumption this variation is not too ``wild''
for the problem at hand),
it is constant in this paper
(except for Remark~\ref{Rem.variable})
and the present complication is that
the underlying manifold is curved.

The paper is organized as follows.
In the next section 
we collect and comment our main results
(Theorems~\ref{Thm.tube.intro} and~\ref{Thm.main.intro})
together with the main ideas behind them.
Section~\ref{Sec.Conv} consists of a number of subsections
devoted to the proof of Theorem~\ref{Thm.main.intro}
concerning the spectral properties of
Schr\"odinger-type operators in a straight tube
of shrinking cross-section
and represents our main technical result.
Section~\ref{Sec.tubes} is devoted to a detailed
definition of a tube and the corresponding Laplacian
(concisely introduced in Section~\ref{Sec.tubes.intro}),
and to an application of Theorem~\ref{Thm.main.intro}
to this situation,
with Theorem~\ref{Thm.tube.intro} as an outcome.
We conclude the paper with Section~\ref{Sec.Strips}
where we discuss an extension of the main result
to strips on surfaces (Theorem~\ref{Thm.strips}).

\setcounter{equation}{0}
\section{Main results and ideas}\label{Sec.Results}
%
Throughout the present paper we use the following notation:
$d \geq 2$ is an integer denoting the dimension;
$I:=(0,L)$ is an open interval of length~$L>0$;
$\cross$ is a bounded open connected subset of~$\Real^{d-1}$
with the centre of mass at the origin
and with the boundary~$\partial\cross$ of class $\Smooth^\infty$;
\begin{equation}\label{def.a}
  a \equiv a(\cross) := \sup_{t\in\omega} |t|
  \,,
\end{equation}
which estimates from above the half of diameter of~$\cross$;
$\cyl:=I\times\cross$ is a $d$-dimensional straight tube
of length~$L$ and cross-section~$\cross$;
and~$\eps$ is a (small) positive parameter.

\subsection{The Laplacian in thin curved tubes}\label{Sec.tubes.intro}
%
We start by some geometric preliminaries
and refer to Section~\ref{Sec.tubes} for more details.

Let $\curve:I\to\Real^d$ be a regular curve
(\ie~an immersion),
which is parametrized by arc length.
We assume that~$\curve$ is uniformly $\Smooth^\infty$-smooth
and that it possesses an appropriate
uniformly $\Smooth^\infty$-smooth Frenet frame $\{e_1,\dots,e_d\}$
(\cf~Assumption~\ref{Ass.Frenet} and Remark~\ref{Rem.Frenet} below).
Then the $i^\mathrm{th}$ curvature~$\kappa_i$ of~$\curve$,
with $i\in\{1,\dots,d-1\}$,
is also uniformly $\Smooth^\infty$-smooth,
and given an (arbitrary) positive constant~$C_\curve$,
we restrict ourselves to the class of curves satisfying
\begin{equation}\label{class}
  \|\kappa_1\|_{\Smooth^2(\overline{I})} \leq C_\curve
  \qquad\mbox{and}\qquad
  \|\kappa_\mu\|_{\Smooth^1(\overline{I})} \leq C_\curve ,
  \,\quad   \forall \mu \in \{2,\dots,d-1\}.
\end{equation}

For any~$\eps>0$,
we introduce the mapping $\tubemap: \Omega \to \Real^d$
by setting
\begin{equation}\label{tube.map.intro}
  \tubemap(s,t) :=
  \curve(s) + \eps \sum_{\mu,\nu=2}^d
  t_\mu \, \rot_{\mu\nu}(s) \, e_\nu(s)
  \,, \qquad
  s \in I, \quad
  t\equiv(t_2,\dots,t_d)\in\cross,
\end{equation}
where $\rot_{\mu\nu}$ are coefficients
of a uniformly $\Smooth^\infty$-smooth
family of rotation matrices in~$\Real^{d-1}$
yet to be specified.
$\tubemap$ is an immersion
provided~$\eps$ is small enough
(namely, it satisfies~(\ref{Ass.basic}) below)
and induces therefore a Riemannian metric~$G$
on the \emph{straight tube}~$\cyl$.
Solving a system of ordinary differential equations
governed by higher curvatures (\cf~(\ref{diff.eq}) below),
we choose the rotations~($\rot_{\mu\nu}$)
in such a special way that the metric is diagonal
(\cf~Section~\ref{Sec.Tang} for more details).
The explicit expression for the metric then reads as follows:
\begin{equation}\label{metric}
  G = \diag\left(h^2,\eps^2,\dots,\eps^2\right)
  \,,
\end{equation}
where the function~$h$ is given by
\begin{equation}\label{Jacobian}
  h(s,t) := 1-\eps \, \kappa_1(s)
  \sum_{\mu=2}^d \rot_{\mu 2}(s) \, t_\mu
  \,,
\end{equation}
and~$\rot_{\mu 2}$ are determined
by the system~(\ref{diff.eq}) below.

\begin{Definition}\label{Def.tubes}
We define the \emph{tube~$\tube$ of cross-section
$
  \eps\cross := \{\eps t\,|\,t\in\cross\}
$
about~$\curve$} to be the Riemannian manifold $(\cyl,G)$,
where the metric~$G$ is given by~(\ref{metric}) with~(\ref{Jacobian}).
\end{Definition}

We refer to Remark~\ref{Rem.rotation} below
for a discussion about the significance
of the special class of rotations, and therefore tubes,
we restrict to.

\begin{Remark}\label{Rem.intersection}
Notice that neither the tube~$\tube$ nor the curve~$\curve$
are required to be embedded in~$\Real^d$.
However, if~$\curve$ is embedded and~$\tubemap$ is injective,
then~$\tubemap$ induces a global diffeomorphism,
$\tube$ is also embedded and the image
$U:=\tubemap(\cyl)$ is an open subset of~$\Real^d$,
which has indeed a geometrical meaning of a non-self-intersecting tube.
Moreover, $\tube$~can be considered as~$U$
expressed in curvilinear ``coordinates" $(s,t)$
via~(\ref{tube.map.intro}).
\end{Remark}

We denote by~$-\Delta_D^\tube$ the \emph{Dirichlet Laplacian}
in the tube $\tube\equiv(\cyl,G)$, \ie,
the self-adjoint operator in the Hilbert space~$\sii(\tube)$
defined as the Friedrichs extension
of the Laplacian on $\Smooth_0^\infty(\tube)$
(\cf~Section~\ref{Sec.Laplacian} for more details).
Of course, in the spirit of Remark~\ref{Rem.intersection},
it is clear that if~$\tubemap$ is injective,
then~$-\Delta_D^\tube$ is nothing else
than the usual Dirichlet Laplacian defined
in the open set~$U$, $-\Delta_D^U$,
and expressed in the ``coordinates" $(s,t)$
via~(\ref{tube.map.intro}).
The spectrum of~$-\Delta_D^\tube$ is purely discrete;
we denote by~$\{\lambda_n\}_{n=1}^\infty$
the set of its eigenvalues
sorted in non-decreasing order and repeated
according to multiplicity,
and by~$\{\Psi_n\}_{n=1}^\infty$ the set
of corresponding eigenfunctions.

As~$\eps \to 0$,
the image of the tube~$U\equiv\tubemap(\cyl)$
collapses in some sense into the reference curve~$\curve$.
However, it turns out that the Dirichlet Laplacian in~$\curve$
(which is in fact the Dirichlet Laplacian in~$I$, $-\Delta_D^I$,
because~$\curve$ is parametrized by arc length)
is not the right operator governing the spectral properties
of~$-\Delta_D^\tube$ in this limit.
Instead, it will be clear in a moment
that the right operator for this purpose is
the one-dimensional Schr\"odinger operator
\begin{equation}\label{op.geom}
  S := -\Delta_D^I + v_0
  \,,
\end{equation}
where the potential function is determined uniquely
by the first curvature of~$\curve$:
\begin{equation}\label{potential.tube}
  v_0 := - \frac{\kappa_1^2}{4}
\end{equation}
(later on, we shall use the notation~$S$
also for other choices of~$v_0$).
The spectrum of~$S$ is purely discrete and simple
for any bounded~$v_0$;
we denote by
$
  \{\mu_n\}_{n=1}^\infty
$
the set of its eigenvalues
sorted in non-decreasing order and repeated according to multiplicity,
and by
$
  \{\phi_n\}_{n=1}^\infty
$
the set of corresponding eigenfunctions.
The Sturm oscillation theorems
imply that each~$\phi_n$ has exactly $n-1$
distinct zeros in~$I$,
which forms therefore the nodal set~$\nodal(\phi_n)$.

Of course, as~$\eps$ decreases to $0$,
the eigenvalues of~$-\Delta_D^\tube$ tend to infinity
because of the Dirichlet boundary conditions,
but it turns out that there is a simple way of
regularizing this singularity.
Namely, the shifted eigenvalues
$
  \lambda_n - \eps^{-2} E_1
$
remain bounded as~$\eps$ goes to $0$, where~$E_1$
is the first eigenvalue of
the Dirichlet Laplacian in~$\cross$, $-\Delta_D^{\cross}$;
we denote by~$\transef_1$ a corresponding eigenfunction.

The connection between the tube~$\tube$
and the operator~$S$ is then given by the following theorem.
\begin{Theorem}\label{Thm.tube.intro}
Given a positive constant~$C_\curve$,
let~$\curve$ be any curve as above satisfying~(\ref{class}),
and let $\tube$ be the tube
of (shrinking) cross-section~$\eps\cross$ about it.
For any integer $N \geq 1$,
there exist positive constants~$\eps_0$ and~$C$
depending on $N,L,C_\curve,\cross$ and~$d$
such that for all $\eps \leq \eps_0$,
\begin{enumerate}
\item[(i)]
the set of eigenvalues $\{\lambda_n\}_{n=1}^\infty$
of~$-\Delta_D^{\tube}$
sorted in non-decreasing order and repeated according to multiplicity
satisfies for any $n\in\{1,\dots,N\}$,
$$
  \big|\lambda_n-(\mu_n+\eps^{-2} E_1)\big|
  \, \leq \, C\,\eps
  \mbox{\,;}
$$
\item[(ii)]
the set of corresponding eigenfunctions $\{\Psi_n\}_{n=1}^N$
of~$-\Delta_D^{\tube}$
can be chosen in such a way that,
for any $n\in\{1,\dots,N\}$,
$$
  \forall (s,t)\in\cyl, \quad
  |\Psi_n(s,t) - \Psi_n^0(s,t)|
  \ \leq \ C\,\eps \ \dist(t,\partial\cross)
  \,,
$$
where
$
  \Psi_n^0 := \eps^{-(d-1)/2} \, h^{-1/2} \, \phi_n\otimes\transef_1
$\,;
\item[(iii)]
with the above choice of eigenfunctions,
we have for any $n\in\{1,\dots,N\}$,
$$
  \forall (s,t)\in\cyl, \quad
  \dist\big(s,\nodal(\phi_n)\big) > C\,\eps
  \ \Longrightarrow \
  \sgn\Psi_n(s,t)
  = \sgn\phi_n(s)
  \,.
$$
\end{enumerate}
\end{Theorem}

As a consequence of~(i), we get that,
as the cross-section of the tube diminishes,
the first~$N$ eigenvalues of~$-\Delta_D^\tube$
converge to the first~$N$ eigenvalues of~$S$
shifted by $\eps^{-2} E_1$.
In particular, the former eigenvalues must also be simple
for all sufficiently small~$\eps$.

Property~(ii) implies that
$\Psi_n$ is well approximated by~$\Psi_n^0$
in the topology of uniformly continuous functions
for all sufficiently small~$\eps$
and that this is also true for ``transverse'' derivatives
on the boundary.

Property~(iii) follows from~(ii)
and the Courant nodal domain theorem
(\cf~Section~\ref{Sec.Nodal} for more details).
In particular, we get that for any $k\in\{1,\dots,n\}$,
the $k^\mathrm{th}$ nodal domain of~$\Psi_n$
projects surjectively to~$\omega$
and converges in some sense to the Cartesian product
of the $k^\mathrm{th}$ nodal domain of~$\phi_n$ with~$\cross$.
At the same time, we get that
the nodal set~$\nodal(\Psi_n)$
projects surjectively to~$\omega$ and converges
to the Cartesian-product $\nodal(\phi_n)\times\cross$.
The Courant nodal domain theorem therefore yields
\begin{Corollary}
Under the hypotheses of Theorem~\ref{Thm.tube.intro},
each $\overline{\nodal(\Psi_n)}$ has a non-empty intersection
with the boundary~$\partial\cyl$
(and this at exactly two points
in the case $d=2$ and $n=2$)
and each~$\Psi_n$ has exactly~$n$ nodal domains.
\end{Corollary}

Consequently, Conjecture~\ref{Conj.simply}
holds for all sufficiently thin tubes
provided the mapping~$\tubemap$ is injective,
since then the above results extend to the Dirichlet Laplacian
in the open set $U \equiv \tubemap(\cyl)$
in view of Remark~\ref{Rem.intersection}.
In particular, if~$u_n$ is the $n^\mathrm{th}$
eigenfunction of~$-\Delta_D^U$
(\eg\ one can choose $u_n:=\Psi_n\circ\tubemap^{-1}$),
the nodal set $\nodal(u_n)$ converges to the image
$
  \tubemap\big(\nodal(\phi_n)\times\cross\big)
$
in the sense that
$$
  \forall x \in U, \quad
  x \in \nodal(u_n)
  \ \Longrightarrow \
  \dist\big(s,\nodal(\phi_n)\big)
  \, \leq \, C\,\eps
  \,,
$$
where~$s$ is determined as the first component
of $\tubemap^{-1}(x)$.

The observation that the operator~$S$
has something to do with the spectral properties
of thin tubes is not new.
Indeed, it is involved in attempts
to justify quantization on submanifolds
(\cf\ \cite{KJ,JK,Tol,daC1,Mitchell} and~\cite{FrHe}
for formal and rigorous treatments, respectively),
in the proofs of the existence of discrete spectrum
in curved quantum wires
(\cf\ \cite{ES,DE,K1,ChDFK} and \cite{ESylwia3}
for different mathematical models, respectively)
and in asymptotics of eigenvalues
of tubular neighbourhoods
(\cf\ \cite{Gray-Karp-Pinsky_1983,Karp-Pinsky_1988}).

Moreover, paper~\cite{DE}
(where the part~(i) of Theorem~\ref{Thm.tube.intro}
is in fact proved for the case $I=\Real$ and $d=2,3$)
inspired us to use $\sii$-perturbation theory
to handle the problem.
However, $\sii$-convergence results are not sufficient
to get an estimate on the location of nodal sets
and we had to develop new ideas
to establish the $\Smooth^0$-
(or even some sort of Lipschitz-)
convergence result included in the part~(ii)
of Theorem~\ref{Thm.tube.intro}.
To our knowledge, the results in the present paper
are the first convergence results of this type
for eigenfunctions and their nodal sets in tubes.

\begin{Remark}[Other tubes]\label{Rem.rotation}
Assume that $\tubemap$ is injective
and recall Remark~\ref{Rem.intersection}.
It should be stressed here that
while the shape of the image~$U$
is not influenced by a special choice of the rotation~$(\rot_{\mu\nu})$
provided~$\cross$ is rotationally symmetric,
this may not be longer true for a general cross-section.
The geometrical meaning of our special choice for rotations~($\rot_{\mu\nu}$)
is that we restrict to \emph{non-twisted} tubes
in the language of~\cite{EKK},
which simplifies the analysis considerably.
It has been noticed recently in~\cite{EKK,BMT}
that other choices for the rotation
may change the spectral picture, too.
Namely, in view of the eigenvalue asymptotics
obtained in~\cite{BMT},
it seems to be reasonable to conjecture
that a version of our Theorem~\ref{Thm.tube.intro}
will still hold for twisted tubes,
provided that~$v_0$ is replaced
by a more complicated potential,
depending also on the higher curvatures of~$\Gamma$
and the geometry of~$\cross$.
\end{Remark}
\begin{Remark}[Tubes of variable cross-section]\label{Rem.variable}
While Theorem~\ref{Thm.tube.intro} applies exclusively
to tubes of uniform cross-section~$\eps\cross$,
it is possible to obtain similar convergence results
for non-uniform cases, too.
Let us fix a tube~$\tube$
of uniform cross-section $\eps\cross$ for which
the results of Theorem~\ref{Thm.tube.intro} hold true.
Consider its deformation~$\tube_\delta$
obtained by replacing~$\eps$ in~(\ref{tube.map.intro})
by $\eps+\delta r(s)$, where~$r$ is a uniformly
$\Smooth^\infty$-smooth function, $\eps$~is fixed
and~$\delta$ plays the role of small parameter now.
Then it is possible to verify by methods of the present paper
that the spectral properties of $-\Delta_D^{\tube_\delta}$
can be approximated by those of~$-\Delta_D^{\tube}$
in the limit $\delta \to 0$,
and obtain an analogy of Theorem~\ref{Thm.tube.intro} in this sense.
However, let us stress that
this problem is much simpler
because the comparison operator is independent
of the perturbation parameter,
while the one treated in the present paper
is singular in~$\eps$.
\end{Remark}
%

\subsection{Schr\"odinger-type operators in a thin straight tube}
\label{Sec.straight.intro}
%
The main idea behind the proof of Theorem~\ref{Thm.tube.intro}
is that the Laplacian $-\Delta_D^\tube$
in a curved tube is unitarily
equivalent to a Schr\"odinger-type operator
in the straight tube~$\cyl$
(\cf~Section~\ref{Sec.Laplacian}).
This observation leads us to consider in Section~\ref{Sec.Conv}
the self-adjoint operator~$T$
in the Hilbert space $\sii(\cyl)$ defined by
\begin{equation}\label{operator}
  T := -\partial_1 a_\eps \partial_1
  + \eps^{-2} (-\Delta'-E_1) + V_\eps
  \,, \qquad
  \Dom(T) := \Hilbert_0^1(\cyl) \cap \Hilbert^2(\cyl)
  \,.
\end{equation}
Here~$-\Delta'$ is the Laplacian
in the variables $(t_2,\dots,t_d)$,
with $(s,t_2,\dots,t_d)\in\cyl$,
$\Hilbert_0^1(\cyl)$ and $\Hilbert^2(\cyl)$
are the usual Sobolev spaces (\cf~\cite{Adams}),
and~$a_\eps$ and~$V_\eps$ are real-valued functions satisfying:
\begin{Assumption}\label{Ass.conv}
There exists a positive constant~$C_{\ref{Ass.conv}}$
(independent of~$\eps$) such that for all
$\eps\in(0,C_{\ref{Ass.conv}}^{-1}]$,
\begin{enumerate}

\item[(i)]
$
  a_\eps, V_\eps \in \Smooth^\infty(\overline{\cyl})
$\,,
\item[(ii)]
$
  \inf_{\Omega} a_\eps \geq C_{\ref{Ass.conv}}^{-1}
$\,,
\item[(iii)]
$
  \|a_\eps-1\|_{\Smooth^1(\overline{\cyl})}
  + \|V_\eps-V_0\|_{\Smooth^0(\overline{\cyl})}
  \leq C_{\ref{Ass.conv}}\,\eps
$\,,
\end{enumerate}
where $V_0:=v_0\otimes\id$
with~$v_0$ being a function from $\Smooth^\infty(\overline{I})$
independent of~$\eps$ and satisfying
\begin{enumerate}
\item[(iv)]
$
  \|v_0\|_{\Smooth^0(\overline{I})} \leq C_{\ref{Ass.conv}}
$\,.
\end{enumerate}
\end{Assumption}
\begin{Remark}
In the case of tubes, $a_\eps$~and~$V_\eps$ will be
expressed in terms of the function~$h$ given by~(\ref{Jacobian})
(\cf~(\ref{Hamiltonian.potential}) and~(\ref{potential.bis}) below)
and~$v_0$ will be given by~(\ref{potential.tube}).
However, we do not restrict ourselves to this particular setting
neither here nor in Section~\ref{Sec.Conv}.
Let us also point out that the subtraction
of~$\eps^{-2}E_1$ in~(\ref{operator})
is a technically useful trick in order to deal
with bounded eigenvalues as $\eps \to 0$
(\cf~Theorem~\ref{Thm.main.intro}.(i) below).
\end{Remark}

The spectrum of~$T$ is purely discrete;
we denote by~$\{\sigma_n\}_{n=1}^\infty$
the set of its eigenvalues
sorted in non-decreasing order and repeated
according to multiplicity,
and by~$\{\psi_n\}_{n=1}^\infty$
the set of corresponding eigenfunctions.
In view of Assumption~\ref{Ass.conv}.(iii),
it is reasonable to expect that, as $\eps \to 0$,
the eigenvalues and eigenfunctions of~$T$
will be approximated by those of the decoupled operator
\begin{equation}\label{decomposition.intro}
  S \otimes \id
  + \id \otimes \eps^{-2}\big(-\Delta_D^{\cross}-E_1\big)
  \qquad\mbox{in}\qquad
  \sii(I)\otimes\sii(\cross)
  \,.
\end{equation}
Here~$S$ is the one-dimensional operator~(\ref{op.geom})
with~$v_0$ being determined only by Assumption~\ref{Ass.conv};
we adopt the notation~$\{\mu_n\}_{n=1}^\infty$
and~$\{\phi_n\}_{n=1}^\infty$
for its set of eigenvalues
(sorted in non-decreasing order and repeated
according to multiplicity)
and corresponding eigenfunctions, respectively;
recall also the notation~$\transef_1$
for the first eigenfunction of $-\Delta_D^\cross$.
It is indeed the case:
\begin{Theorem}\label{Thm.main.intro}
Suppose Assumption~\ref{Ass.conv} holds.
For any integer $N \geq 1$,
there exist positive constants~$\eps_0$ and~$C$
depending on $N,L,C_{\ref{Ass.conv}},\cross$ and~$d$
such that for all $\eps \leq \eps_0$,
\begin{enumerate}
\item[(i)]
the set of eigenvalues $\{\sigma_n\}_{n=1}^\infty$ of~$T$
sorted in non-decreasing order and repeated according to multiplicity
satisfies for any $n\in\{1,\dots,N\}$,
$$
  |\sigma_n-\mu_n|
  \, \leq \, C\,\eps
  \mbox{\,;}
$$
\item[(ii)]
the set of corresponding eigenfunctions $\{\psi_n\}_{n=1}^N$ of~$T$
can be chosen in such a way that for any $n\in\{1,\dots,N\}$,
$$
  \forall (s,t)\in\cyl, \quad
  |\psi_n(s,t) - \psi_n^0(s,t)|
  \ \leq \ C\,\eps \ \dist(t,\partial\cross)
  \,,
$$
where
$
  \psi_n^0 := \phi_n\otimes\transef_1
$\,;
\item[(iii)]
with the above choice of eigenfunctions,
we have for any $n\in\{1,\dots,N\}$,
$$
  \forall (s,t)\in\cyl, \quad
  \dist\big(s,\nodal(\phi_n)\big) > C\,\eps
  \ \Longrightarrow \
  \sgn\psi_n(s,t)
  = \sgn\phi_n(s)
  \,.
$$
\end{enumerate}
\end{Theorem}

Theorem~\ref{Thm.tube.intro} therefore follows
as a consequence of Theorem~\ref{Thm.main.intro}
and the unitary equivalence between
$-\Delta_D^\tube-\eps^{-2}E_1$ and~$T$
with a particular choice of~$a_\eps$ and~$V_\eps$
(\cf~Sections~\ref{Sec.Laplacian} and~\ref{Sec.application}).

The proof of Theorem~\ref{Thm.main.intro}
consists of six steps:
\smallskip \\
1.~Convergence of eigenvalues;
\\
2.~Convergence of eigenfunctions in $\sii(\cyl)$;
\\
3.~Convergence of eigenfunctions in $\Hilbert^2(\cyl)$;
\\
4.~Convergence of eigenfunctions in $\Smooth^0(\cyl)$;
\\
5.~Convergence of transverse derivatives of eigenfunctions
in $\Smooth^0(I\times\partial\cross)$;
\\
6.~Convergence of nodal domains.
\smallskip \\
Step~1, \ie\ part~(i) of Theorem~\ref{Thm.main.intro},
is established in Section~\ref{Sec.evs}
directly through the \emph{minimax principle}.
Step~2 is a consequence of the fact that~$T$
converges to the operator~(\ref{decomposition.intro})
in a generalized sense, which is established
in Section~\ref{Sec.L2} by means of \emph{perturbation theory}.
Step~3 is deduced from Step~2 in Section~\ref{Sec.H2}
by using ideas of \emph{elliptic regularity theory}
in a refined way.
In view of the Sobolev embedding theorem,
Step~3 already implies Step~4 in the case of $d=2$ or~$3$.
For higher dimensions, however,
we have to use a different argument
to deduce Step~4 in Section~\ref{Sec.C0},
namely, the \emph{maximum principle}.
In Section~\ref{Sec.C0},
the latter is also used to establish Step~5,
\ie, more precisely, the part~(ii) of Theorem~\ref{Thm.main.intro}.
Finally, in Section~\ref{Sec.Nodal},
Step~6, \ie\ the part~(iii) of Theorem~\ref{Thm.main.intro},
is deduced from Step~5
by means of the \emph{Courant nodal domain theorem}.

The last reasoning also yields
\begin{Corollary}\label{Corol.main}
Under the hypotheses of Theorem~\ref{Thm.main.intro},
each $\overline{\nodal(\psi_n)}$ has a non-empty intersection
with the boundary~$\partial\cyl$
(and this at exactly two points
in the case $n=2$ and $d=2$)
and each~$\psi_n$ has exactly~$n$ nodal domains.
\end{Corollary}
%

\subsection{Comments on notation}
%
Here we point out some special convention
frequently used throughout the paper.

Since the straight tube~$\cyl$ is of the form $I\times\cross$,
with $I\subset\Real$ and $\cross\subset\Real^{d-1}$,
we consistently split variables
into $(s,t)\in\Real\times\Real^{d-1}$,
where $t\equiv(t_2,\dots,t_d)$.

The symbol $\partial_i^n:=\partial^n/\partial x_i^n$
stands for the partial derivative
of the $n^\mathrm{th}$-order
with respect to the~$i^\mathrm{th}$-variable,
$i\in\{1,\dots,d\}$,
and we use the identification
$
  (x_1,x_2,\dots,x_d)
  \equiv (s,t_2,\dots,t_d)
$.
$\nabla'$ denotes the gradient
in the ``transverse'' variables $(t_2,\dots,t_d)$.

We set $\NatStar:=\Nat\setminus\{0\}$ where $\Nat=\{0,1,\dots\}$.

If~$U$ is an open set, we denote by~$-\Delta_D^U$
the Dirichlet Laplacian in~$U$,
\ie\ the self-adjoint operator associated on~$\sii(U)$
with the quadratic form~$Q_D^U$ defined by
$Q_D^U[\psi]:=\int_U |\nabla\psi|^2$,
$\Dom(Q_D^U):=\Hilbert_0^1(U)$.

We shall usually not distinguish between multiplication operators
and the corresponding generating functions,
between differential operators
and the corresponding differential expressions, \etc.
In fact, all the operators we consider act in the classical sense.

\setcounter{equation}{0}
\section{Schr\"odinger-type operators in shrinking \\
         straight tubes}\label{Sec.Conv}
%
This section is devoted to the study of spectral properties
of the operator~$T$ defined by~(\ref{operator})
in the limit as $\eps$ goes to $0$.
While the cross-section of~$\cyl$ is~$\cross$ (and not $\eps\cross$)
and is thus actually fixed, we refer to~$T$ as a (Schr\"odinger-type)
operator in a shrinking tube, since that is indeed the case
after an obvious change of variables.
Our main goal is Theorem~\ref{Thm.main.intro},
which is established in the following subsections
as described in Section~\ref{Sec.straight.intro}.

\subsection{The Schr\"odinger-type operator}\label{Sec.operator}
%
The operator~$T$ is properly introduced as follows.
Let~$t$ be the sesquilinear form
in the Hilbert space $\Hilbert:=\sii(\cyl)$
defined by
\begin{align*}
  t(\phi,\psi)
  &:= (\partial_1\phi,a_\eps\partial_1\psi)_\Hilbert
  + \eps^{-2} (\nabla'\phi,\nabla'\psi)_\Hilbert
  - \eps^{-2} E_1 (\phi,\psi)_\Hilbert
  + (\phi,V_\eps\psi)_\Hilbert
  \,,
  \\
  \phi,\psi \in \Dom(t)
  &:= \Hilbert_0^1(\cyl)
  \,.
\end{align*}
In view of the properties~(i) and~(ii) from Assumption~\ref{Ass.conv},
the form~$t$ is clearly densely defined, closed, symmetric
and bounded from below for any positive~$\eps$.
Consequently, it follows by the first representation theorem
(\cf~\cite[Sec.~VI.2.1]{Kato}) that there exists
a unique self-adjoint operator~$T$
which is also bounded from below and satisfies
\begin{equation*}
  T\psi = \tau\psi
  \,, \qquad
  \psi\in\Dom(T)
  = \{\psi\in\Hilbert_0^1(\cyl) \,|\ \tau\psi\in\Hilbert\}
  \,,
\end{equation*}
where~$\tau$ is the differential expression in~$\cyl$
defined by
$$
  \tau := -\partial_1 a_\eps \partial_1
  + \eps^{-2} (-\Delta'-E_1) + V_\eps
$$
and the derivatives must be interpreted in the distributional sense.
However, since the boundary~$\partial\cyl$ is sufficiently regular,
it follows by standard methods for regularity of weak solutions
of elliptic equations
(\cf~\cite[Chap.~III]{Lad-Ural})
that indeed
\begin{equation}\label{domain}
  \Dom(T)
  = \Hilbert_2
  := \Hilbert_0^1(\cyl) \cap \Hilbert^2(\cyl)
\end{equation}
and~$T$ acts as~$\tau$ in the classical sense
on functions from $\Dom(T)$.

Hence,
\mbox{$\{T\,|\,\eps\in(0,C_{\ref{Ass.conv}}^{-1}]\}$}
forms a family of elliptic self-adjoint
$\eps$-de\-pen\-dent operators
which is uniformly bounded from below by
$
  - (1+C_{\ref{Ass.conv}})
$
(\cf~(\ref{Poincare}) below for the last statement).
The spectrum of~$T$ is purely discrete
(\cf~\cite[Sec.~8.12]{Gilbarg-Trudinger}),
and recall that we have denoted by~$\{\sigma_n\}_{n=1}^\infty$
the set of its eigenvalues
(sorted in non-decreasing order and repeated
according to multiplicity).
The set of corresponding
eigenfunctions~$\{\psi_n\}_{n=1}^\infty$
can be chosen in such a way that all~$\psi_n$ are real.
Moreover, applying the elliptic regularity theorem repeatedly,
we know that each~$\psi_n\in\Smooth^\infty(\overline{\cyl})$.

\subsection{The comparison operator}
%
In view of Assumption~\ref{Ass.conv}.(iii),
it is reasonable to expect that, as $\eps \to 0$,
$T$~``converges" in a suitable sense
to the self-adjoint operator~$T_0$ defined by
\begin{equation}
  T_0 := -\partial_1^2
  + \eps^{-2} (-\Delta'-E_1) + V_0 \,,
  \qquad
  \Dom(T_0) := \Hilbert_2 \,.
\end{equation}
We use the quotation marks because~$T_0$ itself
forms a family of $\eps$-dependent operators.
$T_0$ is bounded from below by~$-C_{\ref{Ass.conv}}$.
The spectrum of~$T_0$ is purely discrete too
and we denote by~$\{\sigma_n^0\}_{n=1}^\infty$
the set of its eigenvalues sorted in non-decreasing order
and repeated according to multiplicity.

$T_0$ is a good comparison operator since
it is naturally decoupled as~(\ref{decomposition.intro}).
Consequently, we know that (\cf~\cite[Corol.\ of Thm.~VIII.33]{RS1})
\begin{equation}\label{decomposition.sp}
  \sigma(T_0)
  = \{\mu_n\}_{n=1}^\infty + \{\eps^{-2} (E_n-E_1)\}_{n=1}^\infty
  \,,
\end{equation}
where $\{\mu_n\}_{n=1}^\infty$
(respectively~$\{E_n\}_{n=1}^\infty$)
denotes the set of eigenvalues of
the operator~$S$ defined by~(\ref{op.geom})
and with Assumption~\ref{Ass.conv} being the only restriction on~$v_0$
(respectively of $-\Delta_D^{\cross}$),
sorted in non-decreasing order and repeated according to multiplicity.

\subsubsection{The longitudinal operator}\label{Sec.longitudinal}
%
Let~$\phi_n$ be a real eigenfunction of~$S$
corresponding to the eigenvalue~$\mu_n$;
we choose~$\phi_1$ to be positive
and normalize all~$\phi_n$ to~$1$ in $\sii(I)$.
Assumption~\ref{Ass.conv} ensures
that $\phi_n\in\Smooth^\infty(\overline{I})$.
By the Sturm oscillation theorems~\cite[Chap.~X]{Ince},
the spectrum of~$S$ is simple
and~$\phi_n$ has exactly $n+1$ zeros in~$\overline{I}$
which we denote by
\begin{equation}\label{Sturm}
  0 \equiv s_0^0(n) < s_0^1(n) < \dots
  < s_0^{n-1}(n) < s_0^n(n) \equiv L
  \,.
\end{equation}
That is, $\nodal(\phi_n)=\{s_0^k(n)\}_{k=1}^{n-1}$.
We define $I_n:=I\setminus\nodal(\phi_n)$,
\ie~$I$ without the \emph{nodal points} of~$\phi_n$.
We also introduce the open subintervals
$$
  I_n^k := \big(s_0^{k}(n),\,s_0^{k+1}(n)\big)
  \,, \qquad
  k\in\{0,\dots,n-1\}
  \,.
$$

The importance of the following proposition
relies on the fact that
the constant~$c$ is independent of the particular choice of~$v_0$.
\begin{Proposition}\label{Prop.Sturm}
Suppose Assumption~\ref{Ass.conv} holds.
For any $n\in\NatStar$,
there exists a positive constant
$
  c = c(n,L,C_{\ref{Ass.conv}})
$
such that
\begin{enumerate}
\item[(i)]
$
  (n\pi/L)^2 - C_{\ref{Ass.conv}}
  \leq \mu_n \leq
  (n\pi/L)^2 + C_{\ref{Ass.conv}}
$\,;
\item[(ii)]
$
  \|\phi_n\|_{\Smooth^2(\overline{I})} \leq c^{-1}
$\,;
\item[(iii)]
$
  \mu_{n+1}-\mu_n \geq c
$\,;
\item[(iv)]
$
  \forall k\in\{0,\dots,n-1\}, \quad
  s_0^{k+1}(n)-s_0^{k}(n) \geq c
$\,;
\item[(v)]
$
  \forall s\in I, \quad
  |\phi_n(s)| \geq c \, \dist(s,\partial I_n)
$\,.
\end{enumerate}
\end{Proposition}
\begin{proof}
The estimates~(i) follow directly by means of the minimax principle.

The bound~(ii) is deduced in three steps.
Firstly, using the eigenvalue equation for~$S$ and~(i),
we obtain an \emph{a priori} estimate
$
  \|\phi_n\|_{\Hilbert^2(I)} \leq 1+(n\pi/L)^2+2C_{\ref{Ass.conv}}
$.
Secondly,
$
  \|\phi_n\|_{\Smooth^1(\overline{I})}
  \leq C_1 \|\phi_n\|_{\Hilbert^2(I)}
$
with a positive constant
$
  C_1 = C_1(L)
$
coming from the Sobolev embedding theorem
(\cf~\cite[Thm.~5.4]{Adams}).
Finally, using again the eigenvalue equation for~$S$ and~(i),
we can estimate
$
  \|\phi_n''\|_{\Smooth^0(\overline{I})}
  \leq ((n\pi/L)^2 + 2 C_{\ref{Ass.conv}})
  \|\phi_n\|_{\Smooth^0(\overline{I})}
$.
This proves~(ii) with some $c_1=c_1(n,L,C_{\ref{Ass.conv}})>0$.

As for~(iii),
the existence of a positive constant $c_2=c_2(n,L,C_{\ref{Ass.conv}})$
estimating from below the eigenvalue gaps
uniformly in the class of potentials
bounded by a constant~$C_{\ref{Ass.conv}}$
follows from~\cite{Kirsch-Simon_1985}.
(Alternatively, one can use the Rellich-Kondrachov
embedding theorem in the spirit of
\cite[proof of Thm.~II.1]{Ashbaugh-Harrell-Svirsky_1991},
where~(iii) is proved for $n=1$;
the generalization to higher eigenvalues is straightforward.)

To prove~(iv), notice that
for any given $n\in\NatStar$,
$\phi_n$ satisfies the eigenvalue problem
\begin{equation}\label{subintervals}
  \left\{
  \begin{aligned}
  -\phi_n''+v_0\phi_n&=\mu_n\phi_n
  \quad &\mbox{in}& \quad I_n^k
  \,,
  \\
  \phi_n &= 0
  \quad &\mbox{on}& \quad \partial I_n^k
  \,,
  \end{aligned}
  \right.
\end{equation}
for every $k\in\{0,\dots,n-1\}$.
Furthermore, $\phi_n$~does not change sign in~$I_n^k$.
Combining~(\ref{subintervals}) with~(i),
we arrive at
$$
  n^2\pi^2/L^2 + C_{\ref{Ass.conv}}
  \geq \mu_n \geq
  \pi^2/\big(s_0^{k+1}(n)-s_0^{k}(n)\big)^2 - C_{\ref{Ass.conv}}
  \,,
$$
which gives~(iv)
with a positive constant~$c_3=c_3(n,L,C_{\ref{Ass.conv}})$.

Finally, given $n\in\NatStar$, fix $k\in\{0,\dots,n-1\}$
and let~$s_\mathrm{max}\in I_n^k$
be such that $\max_{I_n^k}|\phi_n|=|\phi_n(s_\mathrm{max})|$.
Combining~(\ref{subintervals}) with~(i)
and using obvious estimates,
we arrive at
\begin{align*}
  \pi^2/\big(s_0^{k+1}(n)-s_0^{k}(n)\big)^2
  & \leq \int_{s_0^k(n)}^{s_0^{k+1}(n)} \phi_n'(s)^2 ds
  = \int_{s_0^k(n)}^{s_0^{k+1}(n)} [\mu_n-v_0(s)]\phi_n(s)^2 ds
  \\
  &\leq (n^2\pi^2/L^2 + 2 C_{\ref{Ass.conv}}) \,
  \big(s_0^{k+1}(n)-s_0^{k}(n)\big)
  \, \phi_n(s_\mathrm{max})^2
  \,.
\end{align*}
Since $s_0^{k+1}(n)-s_0^{k}(n) \leq L$,
it follows that
$
  |\phi_n(s_\mathrm{max})| \geq c_4
$
with a positive constant $c_4=c_4(n,L,C_{\ref{Ass.conv}})$.
Let~$w$ be the solution to
\begin{equation*}
  \left\{
  \begin{aligned}
  A w:=-w''+C_{\ref{Ass.conv}}w &= 0
  \quad &\mbox{in}& \quad (s_0^{k}(n),\,s_\mathrm{max})
  \,,
  \\
  w &= 0
  \quad &\mbox{at}& \quad s_0^{k}(n)
  \,,
  \\
  w &= c_4
  \quad &\mbox{at}& \quad s_\mathrm{max}
  \,.
  \end{aligned}
  \right.
\end{equation*}
Since $Aw \leq A|\phi_n|$ in $(s_0^{k}(n),\,s_\mathrm{max})$
and $w \leq |\phi_n|$ at the boundary points,
the maximum principle
(\cf~\cite[Thm.~3.3]{Gilbarg-Trudinger}) yields
$w \leq |\phi_n|$ in $(s_0^{k}(n),\,s_\mathrm{max})$.
Using the explicit form of~$w$, we estimate
$$
  w(s) \equiv \frac{c_4 \sinh[\sqrt{C_{\ref{Ass.conv}}}(s-s_0^{k}(n))]}
  {\sinh[\sqrt{C_{\ref{Ass.conv}}}(s_\mathrm{max}-s_0^{k}(n))]}
  \geq \frac{c_4\sqrt{C_{\ref{Ass.conv}}}\,(s-s_0^{k}(n))}
  {\sinh[\sqrt{C_{\ref{Ass.conv}}}L]}
  \geq \frac{c_4\sqrt{C_{\ref{Ass.conv}}}\,\dist(s,\partial I_n)}
  {\sinh[\sqrt{C_{\ref{Ass.conv}}}L]}
  \,,
$$
and therefore $|\phi_n| \geq c_5 \dist(\cdot,\partial I_n)$
in $(s_0^{k}(n),\,s_\mathrm{max})$
with a positive constant~$c_5=c_5(n,L,C_{\ref{Ass.conv}})$.
A similar comparison argument on $(s_\mathrm{max},s_0^{k+1}(n))$
yields the same lower bound, with the same constant~$c_5$,
and~(v) is proved.

The constant $c$ is eventually chosen
as the smallest one among $c_1,\dots,c_5$.
\end{proof}
%

\subsubsection{The transversal operator}
%
Let~$\transef_n$ denote a real eigenfunction of~$-\Delta_D^\cross$
corresponding to~$E_n$; we choose~$\transef_1$ to be positive
and normalize all~$\transef_n$ to~$1$ in $\sii(\cross)$.
Since we assume that the boundary~$\partial\cross$
is of class~$\Smooth^\infty$,
elliptic regularity theory yields
that the eigenfunctions
belong to $\Smooth^\infty(\overline{\cross})$
(they are even analytic in the interior of~$\cross$).
Moreover, since~$\transef_1$ is a positive eigenfunction
of a Dirichlet problem, we have
\begin{Proposition}\label{Prop.trans}
There exists a positive constant
$
  c = c(\cross)
$
such that
$$
  \forall t\in\cross, \qquad
  \transef_1(t) \, \geq \, c \,\, \dist(t,\partial\cross)
  \,.
$$
\end{Proposition}

Next, as a consequence of the variational definition of~$E_1$,
we have the following Poincar\'e inequality for~$\cross$:

\begin{equation}\label{Poincare}
  \forall\psi\in\Hilbert_0^1(\cross), \qquad
  \|\nabla\psi\|_{\sii(\cross)}^2
  \,\geq\, E_1 \, \|\psi\|_{\sii(\cross)}^2 \,.
\end{equation}
This inequality will play a major role in what follows.

\subsubsection{Spectrum of the comparison operator}
%
\begin{Proposition}\label{Prop.simple}
Suppose Assumption~\ref{Ass.conv} holds.
One has $\sigma_1^0 = \mu_1$.
Moreover, for any integer $N \geq 2$, there exists
$
  \eps_0
  = \eps_0(N,L,C_{\ref{Ass.conv}},\cross) > 0
$
such that for all $\eps \leq \eps_0$,
\begin{equation*}
  \forall n\in\{1,\dots,N\}, \qquad
  \sigma_n^0 = \mu_n
  \,.
\end{equation*}
\end{Proposition}
\begin{proof}
Recall~(\ref{decomposition.sp}).
The assertion for~$n=1$ is obvious.
Let~$n \geq 2$ and assume by induction
that $\sigma_{n-1}^0 = \mu_{n-1}$.
Then
$
  \sigma_{n}^0 = \min\left\{
  \mu_{n}, \mu_{n-1} + \eps^{-2}(E_2 - E_1)
  \right\}
$
and the assertion of Proposition follows at once.
It is only important to notice that
$E_2 - E_1$ depends on~$\cross$,
and $\mu_{n}-\mu_{n-1}$
can be estimated by means of Proposition~\ref{Prop.Sturm}.(i).
\end{proof}

Let~$\psi_n^0$ denote a real eigenfunction of~$T_0$
corresponding to~$\sigma_n^0$.
In view of Proposition~\ref{Prop.simple}
and the fact that~$T_0$ is decoupled as~(\ref{decomposition.intro}),
there is a natural choice for~$\psi_n^0$
if~$\eps$ is small enough. Namely, we always choose
\begin{equation}\label{ef}
  \psi_n^0 := \phi_n \otimes \transef_1
\end{equation}
provided the conclusion of Proposition~\ref{Prop.simple}
holds for the first~$n$ eigenvalues.
$\psi_n^0$~is then normalized to~$1$ in~$\Hilbert$.
As a consequence of~(\ref{ef}),
Proposition~\ref{Prop.Sturm}.(v) and Proposition~\ref{Prop.trans},
we get
\begin{Proposition}\label{Prop.nodal}
Suppose Assumption~\ref{Ass.conv} holds.
For any  $N\in\NatStar$,
there exist positive constants
$
  \eps_0
  = \eps_0(N,L,C_{\ref{Ass.conv}},\cross)
$
and
$
  c
  = c(N,L,C_{\ref{Ass.conv}},\cross)
$
such that for all $\eps \leq \eps_0$
and any $n\in\{1,\dots,N\}$,
\begin{enumerate}
\item[(i)]
$
  \nodal(\psi_n^0) = \nodal(\phi_n) \times \cross
$\,,
\item[(ii)]
$
  \forall (s,t)\in\cyl, \quad
  |\psi_n^0(s,t)|
  \,\geq\, c \,\,
  \dist(s,\partial I_n) \,\, \dist(t,\partial\cross)
$\,.
\end{enumerate}
\end{Proposition}
%

\subsection{Convergence of eigenvalues}\label{Sec.evs}
%
The purpose of this subsection is to show
that~$T$ converges to~$T_0$ in the sense of their spectra.
\begin{Theorem}\label{Thm.evs}
Suppose Assumption~\ref{Ass.conv} holds.
For any~$N\in\NatStar$,
there exist positive constants
$
  \eps_0
  = \eps_0 (N,L,C_{\ref{Ass.conv}},\cross)
$
and
$
  C
  = C(N,L,C_{\ref{Ass.conv}})
$
such that for all $\eps \leq \eps_0$,
$$
  \forall n\in\{1,\dots,N\}, \qquad
  |\sigma_n-\sigma^0_n|
  \ \leq \ C \,\eps
  \,.
$$
\end{Theorem}
\begin{proof}
Using Assumption~\ref{Ass.conv},
we estimate
$
  T^- \leq T \leq T^+
$
in the form sense,
where~$T^\pm$ are defined by $\Dom(T^\pm):=\Hilbert_2$ and
$$
  T^\pm := (1 \pm C_{\ref{Ass.conv}} \eps) (-\partial_1^2+V_0)
  + \eps^{-2} (-\Delta'-E_1)
  \pm C_{\ref{Ass.conv}}(1+C_{\ref{Ass.conv}}) \eps
  \,.
$$
Assuming that~$\eps<C_{\ref{Ass.conv}}^{-1}$,
so that~$T^-$ is bounded from below,
the minimax principle gives
$
  \sigma_n^-\leq\sigma_n\leq\sigma_n^+
$
for all $n\in\NatStar$,
where~$\{\sigma_n^\pm\}_{n=1}^\infty$
denotes the set of eigenvalues of~$T^\pm$
sorted in non-decreasing order and repeated
according to multiplicity.
However, using a similar argument
to that leading to Proposition~\ref{Prop.simple},
we know that if~$\eps$ is so small that the inequality
$
  (1+C_{\ref{Ass.conv}} \eps) \, \eps^2
  \leq (E_2-E_1)/(\mu_n-\mu_1)
$
holds true,
where $E_2 - E_1$ depends on~$\cross$ and $\mu_{n}-\mu_{n-1}$
can be estimated by means of Proposition~\ref{Prop.Sturm}.(i),
then
$$
  \sigma_n^\pm
  = (1 \pm C_{\ref{Ass.conv}} \eps) \mu_n
  \pm C_{\ref{Ass.conv}}(1+C_{\ref{Ass.conv}}) \eps
  = \sigma_n^0
  \pm C_{\ref{Ass.conv}}(\mu_n+1+C_{\ref{Ass.conv}}) \eps
  \,.
$$
This proves the claim, since~$\mu_n$
can be estimated by Proposition~\ref{Prop.Sturm}.(i).
\end{proof}
\begin{Remark}
The $\Smooth^1$-norm of $a_\eps-1$
in Assumption~\ref{Ass.conv}.(iii) could be weakened
to the mere $\Smooth^0$-norm in order to prove Theorem~\ref{Thm.evs}
by the above method.
\end{Remark}

The first part of Theorem~\ref{Thm.main.intro}
follows immediately as a consequence
of Theorem~\ref{Thm.evs} and Proposition~\ref{Prop.simple}.

\subsection{$L^2$-convergence of eigenfunctions}\label{Sec.L2}
%
While we have proved Theorem~\ref{Thm.evs} using just the minimax principle,
we shall have to use a stronger technique
to establish the convergence of eigenfunctions in~$\Hilbert$,
namely, perturbation theory.

Let us start by introducing some notation.
Given two Hilbert spaces~$X$ and~$X'$,
we denote by $\Bounded(X,X')$
the set of bounded operators from~$X$ to~$X'$;
we also denote $\Bounded(X):=\Bounded(X,X)$.
For every complex number~$z$ in the resolvent set of~$T$
[respectively~$T_0$],
we introduce $R(z):=(T-z)^{-1}\in\Bounded(\Hilbert)$
[respectively $R_0(z):=(T_0-z)^{-1}\in\Bounded(\Hilbert)$].
The space~$\Hilbert_2$ is introduced in~(\ref{domain})
and we equip it with the usual $\Hilbert^2(\cyl)$-norm.

We shall need two technical lemmata.
\begin{Lemma}\label{Lem.noc}
Suppose Assumption~\ref{Ass.conv} holds.
One has
$$
  \|T-T_0\|_{\Bounded(\Hilbert_2,\Hilbert)}
  \leq 2 \, C_{\ref{Ass.conv}} \, \eps
  \,.
$$
\end{Lemma}
\begin{proof}
Using Assumption~\ref{Ass.conv}, one has the obvious estimates
\begin{multline*}
  \|(T-T_0)\psi\|_{\Hilbert}^2
  \equiv \|(1-a_\eps)\partial_1^2\psi-(\partial_1 a_\eps)\partial_1\psi
  +(V_\eps-V_0)\psi\|_\Hilbert^2
  \\
  \leq 4\,C_{\ref{Ass.conv}}^2 \, \eps^2 \,
  \big(
  \|\partial_1^2\psi\|_\Hilbert^2
  +\|\partial_1\psi\|_\Hilbert^2
  +\|\psi\|_\Hilbert^2
  \big)
  \leq 4\, C_{\ref{Ass.conv}}^2 \, \eps^2 \, \|\psi\|_{\Hilbert_2}^2
\end{multline*}
for any $\psi\in\Hilbert_2$.
\end{proof}
\begin{Lemma}\label{Lem.resolvent}
Suppose Assumption~\ref{Ass.conv} holds.
There exists a negative number
$
  z_0
  = z_0(C_{\ref{Ass.conv}},\cross)
$
such that for all $\eps \leq 1$,
$$
  \forall z \leq z_0, \qquad
  \|R_0(z)\|_{\Bounded(\Hilbert,\Hilbert_2)} \leq 2
  \,.
$$
\end{Lemma}
\begin{proof}
For any $z\in\Real\setminus\sigma(T_0)$
and $f\in\Smooth_0^\infty(\cyl)$,
let~$\psi$ be the (unique) solution to
\begin{equation*}
  \left\{
  \begin{aligned}
  (T_0 - z)\psi &= f
  &\mbox{in}& \quad \cyl
  \,,
  \\
  \psi &= 0
  &\mbox{on}& \quad \partial\cyl
  \,.
  \end{aligned}
  \right.
\end{equation*}
Assumption~\ref{Ass.conv} ensures that
$\psi\in\Smooth^\infty(\overline{\cyl})$.
Restricting ourselves to $z \leq - C_{\ref{Ass.conv}}$
so that $V_0-z \geq 0$,
and integrating by parts, we have
\begin{align*}
  \|f\|_\Hilbert^2
  \ = \ & \|(T_0 - z)\psi\|_\Hilbert^2
  \\
  \ = \ & \|\partial_1^2\psi\|_\Hilbert^2
  + \|(V_0-z)\psi\|_\Hilbert^2
  - 2 z \|\partial_1\psi\|_\Hilbert^2
  - 2 \Re (\partial_1^2\psi,V_0\psi)_\Hilbert
  \\
  & + 2\eps^{-2} \big(
  \|\nabla'\partial_1\psi\|_\Hilbert^2
  - E_1 \|\partial_1\psi\|_\Hilbert^2
  \big)
  \\
  & + 2\eps^{-2} \big(
  \|\nabla'(V_0-z)^{1/2}\psi\|_\Hilbert^2
  - E_1 \|(V_0-z)^{1/2}\psi\|_\Hilbert^2
  \big)
  \\
  & + \eps^{-4} \big(
  \|\Delta'\psi\|_\Hilbert^2
  - 2 E_1 \|\nabla'\psi\|_\Hilbert^2
  + E_1^2 \|\psi\|_\Hilbert^2
  \big)
  \,.
\end{align*}
It is important to notice that the last three lines are non-negative.
This is clear for the last one since the expression in brackets
is just $\|(-\Delta'-E_1)\psi\|_\Hilbert^2$
after a factorization and an integration by parts,
while the precedent two are non-negative due to~(\ref{Poincare}).
We also note that
$
  2 |\Re (\partial_1^2\psi,V_0\psi)_\Hilbert|
  \leq \mbox{$\demi$} \|\partial_1^2\psi\|_\Hilbert^2
  + 2 \|V_0\psi\|_\Hilbert^2
$
by a consecutive usage of the Schwarz and Cauchy inequalities.
Consequently,
restricting ourselves to $\eps \leq 1$,
rearranging the norms,
recalling that $\nabla'V_0=0$
and using obvious estimates,
we can write
\begin{align*}
  \|f\|_\Hilbert^2
  \ \geq \ & \mbox{$\demi$} \|\partial_1^2\psi\|_\Hilbert^2
  + 2 \|\nabla'\partial_1\psi\|_\Hilbert^2
  + \|\Delta'\psi\|_\Hilbert^2
  \\
  & + 2 (-z-E_1) \|\partial_1\psi\|_\Hilbert^2
  + 2 (-z-\|V_0\|_{\sinf(\cyl)}-E_1) \|\nabla'\psi\|_\Hilbert^2
  \\
  & + \big[
  (-z-\|V_0\|_{\sinf(\cyl)}+E_1)^2
  - 2 \|V_0\|_{\sinf(\cyl)}^2
  \big]
  \|\psi\|_\Hilbert^2
  \,.
\end{align*}
Choosing~$z$ sufficiently large negative,
namely
$
  z \leq z_0 := -(1 + 3 C_{\ref{Ass.conv}} + E_1)
$,
we therefore conclude by
\begin{equation*}
  \|f\|_\Hilbert^2
  \ \geq \
  \mbox{$\demi$} \|\Delta\psi\|_\Hilbert^2
  + \|\nabla\psi\|_\Hilbert^2 + \|\psi\|_\Hilbert^2
  \ \geq \
  \mbox{$\demi$} \|\psi\|_{\Hilbert_2}^2
  \,.
\end{equation*}
Then the assertion of Lemma
follows by the density of $\Smooth_0^\infty(\cyl)$
in~$\Hilbert$.
\end{proof}

Now we are ready to prove the main result of this subsection.
\begin{Theorem}\label{Thm.nrc}
Suppose Assumption~\ref{Ass.conv} holds.
Given $n\in\NatStar$, let~$\mathcal{C}$ be the circle in~$\Com$
centred in~$\mu_n$ and of radius~$r$ defined by
$$
  4r :=
  \begin{cases}
    \min\big\{
    c(n-1,L,C_{\ref{Ass.conv}}),
    c(n,L,C_{\ref{Ass.conv}})
    \big\}
    & \mbox{if} \quad n \geq 2 \,,
    \\
    c(1,L,C_{\ref{Ass.conv}})
    & \mbox{if} \quad n = 1 \,,
  \end{cases}
$$
where~$c$ is the positive constant determined
by Proposition~\ref{Prop.Sturm}.
Then there exist positive constants
$
  \eps_0
  = \eps_0(n,L,C_{\ref{Ass.conv}},\cross)
$
and
$
  C = C(n,L,C_{\ref{Ass.conv}})
$
such that for all $\eps \leq \eps_0$,
one has
$
  \mathcal{C}\cap[\sigma(T_0)\cup\sigma(T)]
  = \varnothing
$
and
$$
  \forall z\in\mathcal{C}, \qquad
  \|R(z)-R_0(z)\|_{\Bounded(\Hilbert)}
  \ \leq \ C\,\eps
  \,.
$$
\end{Theorem}
\begin{proof}
Fix $n\in\NatStar$.
Let~$\eps$ be sufficiently small so
that the conclusion of Proposition~\ref{Prop.simple}
holds true for~$n+1$ eigenvalues of~$T_0$.
In particular,
$
  \dist(\sigma(T_0),\mathcal{C}) = r
$.
We also assume that $\eps \leq 1$,
so that the conclusion of Lemma~\ref{Lem.resolvent} holds true
with some $z_0<0$.
Then the first resolvent identity~\cite[Thm.~5.13.(a)]{Weidmann},
the embedding $\Hilbert_2\hookrightarrow\Hilbert$
and obvious estimates yield for any $z\in\mathcal{C}$,
\begin{align*}
  \|R_0(z)\|_{\Bounded(\Hilbert,\Hilbert_2)}
  &\leq \|R_0(z_0)\|_{\Bounded(\Hilbert,\Hilbert_2)} \,
  \big(
  1+(|z|+|z_0|) \, \|R_0(z)\|_{\Bounded(\Hilbert)}
  \big)
  \\
  &\leq
  2 \big(
  1+(\mu_n+r+|z_0|) \, r^{-1}
  \big)
  \leq C_1
\end{align*}
where $C_1=C_1(n,L,C_{\ref{Ass.conv}})$
comes from an upper bound to~$\mu_n$
due to Proposition~\ref{Prop.Sturm}.(i).
Let us now assume in addition that~$\eps$ is sufficiently small so
that~$\sigma_n$ lies inside
the circle concentric with~$\mathcal{C}$
but with half the radius and other eigenvalues of~$T$
lie outside the circle concentric with~$\mathcal{C}$
but with twice the radius
(this will be true provided
$
  \eps \leq \eps_0(n+1,L,C_{\ref{Ass.conv}},\cross)
$
and
$
  C(n+1,L,C_{\ref{Ass.conv}})\,\eps < r/2
$,
where~$\eps_0$ and~$C$ are determined
by Theorem~\ref{Thm.evs}).
In particular,
$
  \dist(\sigma(T),\mathcal{C}) \geq r/2
$.
Then the second resolvent identity \cite[Thm.~5.13.(c)]{Weidmann}
and Lemma~\ref{Lem.noc} yield
\begin{align*}
  \|R(z)-R_0(z)\|_{\Bounded(\Hilbert)}
  &\leq \|R(z)\|_{\Bounded(\Hilbert)} \,
  \|T-T_0\|_{\Bounded(\Hilbert_2,\Hilbert)} \,
  \|R_0(z)\|_{\Bounded(\Hilbert,\Hilbert_2)}
  \\
  &\leq 2 \, r^{-1} \ 2 \, C_{\ref{Ass.conv}} \, \eps \ C_1
  \,,
\end{align*}
which concludes the proof.
\end{proof}

As a consequence of Theorem~\ref{Thm.nrc},
we get that, as $\eps \to 0$,
$T$~converges to~$T_0$ in the \emph{generalized sense}
\cite[Sec.~IV.2.6]{Kato}
(also referred to as the
\emph{norm resolvent sense}~\cite[Sec.~2.6]{Davies}).
This implies (\cf~\cite[Sec.~IV.3.5]{Kato}) the continuity
of eigenvalues of~$T$ at~$\eps=0$ (\cf~also our Theorem~\ref{Thm.evs})
and of the corresponding spectral projections
$$
  P_n := -\frac{1}{2 \pi i} \int_\mathcal{C} R(z) \, dz
  \,,
$$
where~$\mathcal{C}$ is determined in Theorem~\ref{Thm.nrc}
for a given $n\in\NatStar$ and $\eps\leq\eps_0$.
Furthermore, one gets the continuity of eigenfunctions
if they are normalized suitably;
we choose
\begin{equation}\label{efs}
  \psi_n := \Re\left(P_n\psi_n^0\right)
\end{equation}
and stress that~$\psi_n$ is not normalized to~1 in~$\Hilbert$
with this choice.
\begin{Corollary}\label{Thm.efs}
Suppose Assumption~\ref{Ass.conv} holds.
For any~$N\in\NatStar$, there exist
positive constants
$
  \eps_0
  = \eps_0(N,L,C_{\ref{Ass.conv}},\cross)
$
and
$
  C = C(N,L,C_{\ref{Ass.conv}})
$
such that for all $\eps\leq\eps_0$,
$$
  \forall n\in\{1,\dots,N\}, \qquad
  \|\psi_n-\psi^0_n\|_{\Hilbert}
  \ \leq \ C\,\eps
  \,.
$$
\end{Corollary}
%

\subsection{$\Hilbert^2$-convergence of eigenfunctions}\label{Sec.H2}
%
Since~$\psi_n$ and~$\psi_n^0$ are eigenfunctions
of elliptic operators, it is possible to deduce
from Corollary~\ref{Thm.efs} a stronger convergence result.
To do so, let us point out several facts.
Fix $n\in\NatStar$.
Combining the eigenvalue equations for~$T$ and~$T_0$,
we verify that the difference
\begin{equation}\label{psi}
  \psi:=\psi_n-\psi_n^0
\end{equation}
satisfies the Dirichlet problem
\begin{equation}\label{identity}
  \left\{
  \begin{aligned}
  -\partial_1\!\left(a_\eps\,\partial_1 \psi\right)
  + \eps^{-2} \left(-\Delta'-E_1\right)\psi
  + b_\eps\psi
  &= f_\eps
  \quad &\mbox{in}& \quad \cyl \,,
  \\
  \psi &= 0
  \quad &\mbox{on}& \quad \partial\cyl \,,
  \end{aligned}
  \right.
\end{equation}
where
$$
  b_\eps :=V_\eps-\sigma_n \,, \qquad
  f_\eps := (a_\eps-1) \partial_1^2\psi_n^0
  + (\partial_1 a_\eps) \partial_1\psi_n^0
  + (\sigma_n-\sigma_n^0-V_\eps+V_0) \psi_n^0 \,.
$$
The functions~$b_\eps$ and~$f_\eps$
belong to $\Smooth^\infty(\overline{\Omega})$
due to Assumption~\ref{Ass.conv}.
Moreover, using in addition~(\ref{ef}),
Proposition~\ref{Prop.Sturm} and Theorem~\ref{Thm.evs},
it is easy to check that
\begin{equation}\label{Ass.conv'}
  \|b_\eps\|_{\Smooth^0(\overline{\Omega})}
  \leq C_{\ref{Ass.conv}}'
  \qquad\mbox{and}\qquad
  \|f_\eps\|_{\Smooth^0(\overline{\Omega})}
  \leq  C_{\ref{Ass.conv}}' \,\eps
\end{equation}
with some positive
$
  C_{\ref{Ass.conv}}'=C_{\ref{Ass.conv}}'(n,L,C_{\ref{Ass.conv}},\cross)
$.
We also recall that $\psi\in\Smooth^\infty(\overline{\Omega})$
due to Assumption~\ref{Ass.conv}.
Now we are ready to prove
\begin{Theorem}\label{Thm.H2c}
Suppose Assumption~\ref{Ass.conv} holds.
For any~$N\in\NatStar$, there exist
positive constants
$
  \eps_0
  = \eps_0(N,L,C_{\ref{Ass.conv}},\cross)
$
and
$
  C = C(N,L,C_{\ref{Ass.conv}},\cross)
$
such that for all $\eps \leq \eps_0$,
$$
  \forall n\in\{1,\dots,N\}, \qquad
  \|\psi_n-\psi^0_n\|_{\Hilbert^2(\cyl)}
  \ \leq \ C\,\eps
  \,.
$$
\end{Theorem}
\begin{proof}
Fix $n\in\NatStar$ and assume that~$\eps$
is sufficiently small so that the conclusion of Corollary~\ref{Thm.efs}
holds true with a positive constant~$C_1$.

Multiplying the first equation of~(\ref{identity}) by~$\psi$,
integrating by parts in~$\Omega$
and using obvious estimates,
we arrive at
\begin{multline*}
  \big\|a_\eps^{1/2}\partial_1\psi\big\|_\Hilbert^2
  + \eps^{-2} \big(\|\nabla'\psi\|_\Hilbert^2-E_1\|\psi\|_\Hilbert^2\big)
  \\
  = -(\psi,b_\eps\psi)_\Hilbert + (\psi,f_\eps)_\Hilbert
  \leq \|b_\eps\|_{\Smooth^0(\overline{\Omega})} \|\psi\|_\Hilbert^2
  + \|\psi\|_\Hilbert^2 + \|f_\eps\|_\Hilbert^2
  \leq C_2^2 \, \eps^2
  \,,
\end{multline*}
where $C_2=C_2(n,L,C_{\ref{Ass.conv}},\cross)$
is a positive constant determined by~$C_{\ref{Ass.conv}}'$,
$C_1$ and the volume of~$\cyl$.
Since
$
  \big\|a_\eps^{1/2}\partial_1\psi\big\|_\Hilbert^2
  \geq C_{\ref{Ass.conv}}^{-1} \|\partial_1\psi\|_\Hilbert^2
$
by Assumption~\ref{Ass.conv}.(ii)
and
$
  \|\nabla'\psi\|_\Hilbert^2 - E_1\|\psi\|_\Hilbert^2
$
is non-negative by~(\ref{Poincare}), we get that
\begin{equation}\label{H1}
  \|\partial_1\psi\|_\Hilbert \leq C_3\,\eps
  \,, \qquad
  \|\nabla'\psi\|_\Hilbert \leq C_3\,\eps
  \,,
\end{equation}
where $C_3=C_3(n,L,C_{\ref{Ass.conv}},\cross)$
is a positive constant determined
by~$C_2$, $C_{\ref{Ass.conv}}$, $C_1$ and~$E_1$.
This proves an $\Hilbert^1$-convergence of~$\psi$.

Now we rewrite the first equation of~(\ref{identity}) as
\begin{equation}\label{identity.bis}
  -\partial_1^2\psi
  + \eps^{-2} \left(-\Delta'-E_1\right)\psi
  = (a_\eps-1) \partial_1^2\psi + f_\eps'
  \,,
\end{equation}
where $f_\eps':=f_\eps -b_\eps\psi + (\partial_1 a_\eps) \partial_1\psi$.
Taking the norm of both sides of~(\ref{identity.bis})
and using obvious estimates,
we arrive at
\begin{eqnarray*}
\lefteqn{\|\partial_1^2\psi\|_\Hilbert^2
  - 2 \eps^{-2} \big(\partial_1^2\psi,(-\Delta'-E_1)\psi\big)_\Hilbert
  + \eps^{-4} \big\|(-\Delta'-E_1)\psi\big\|_\Hilbert^2}
  \rule{30ex}{0ex}
  \\
  = \big\|(a_\eps-1) \partial_1^2\psi + f_\eps'\big\|_\Hilbert^2
  &\leq& 2 \|a_\eps-1\|_{\Smooth^0(\overline{\Omega})}^2
  \|\partial_1^2\psi\|_\Hilbert^2
  + 2 \|f_\eps'\|_\Hilbert^2
  \\
  &\leq& C_4^2\,\eps^2 \big(\|\partial_1^2\psi\|_\Hilbert^2 + 1 \big)
  \,,
\end{eqnarray*}
where $C_4=C_4(n,L,C_{\ref{Ass.conv}},\cross)$
is a positive constant determined by $C_{\ref{Ass.conv}}$,
$C_{\ref{Ass.conv}}'$, $C_1$, $C_3$
and the volume of~$\cyl$.
Since an integration by parts yields
\begin{align*}
  -\big(\partial_1^2\psi,(-\Delta'-E_1)\psi\big)_\Hilbert
  &= \|\nabla'\partial_1\psi\|_\Hilbert^2
  - E_1 \|\partial_1\psi\|_\Hilbert^2 \,,
  \\
  \big\|(-\Delta'-E_1)\psi\big\|_\Hilbert^2
  &= \|\Delta'\psi\|_\Hilbert^2
  - 2 E_1 \|\nabla'\psi\|_\Hilbert^2
  + E_1^2 \|\psi\|_\Hilbert^2 \,,
\end{align*}
where the expression on the right hand side of the first line
is non-negative by~(\ref{Poincare}),
we conclude that if $C_4\,\eps<1$ then
\begin{equation}\label{H2}
  \|\partial_1^2\psi\|_\Hilbert \leq C_5\,\eps
  \,, \qquad
  \|\nabla'\partial_1\psi\|_\Hilbert \leq C_5\,\eps
  \,, \qquad
  \|\Delta'\psi\|_\Hilbert \leq C_5\,\eps
  \,,
\end{equation}
where $C_5=C_5(n,L,C_{\ref{Ass.conv}},\cross)$
is a positive constant determined by~$C_4$, $C_3$,
$C_1$ and~$E_1$.
Summing up, Corollary~\ref{Thm.efs}, (\ref{H1}) and~(\ref{H2})
establish the assertion of Theorem.
\end{proof}
%

\subsection{$\Smooth^0$-convergence of eigenfunctions}\label{Sec.C0}
%
To show that the convergence result of Theorem~\ref{Thm.H2c}
holds actually in the topology of $\Smooth^0(\overline{\cyl})$
in all dimensions,
we shall use the fact that~$\psi$
is a classical solution of~(\ref{identity})
satisfying the maximum principle.
In particular, our method is based on:
\begin{Proposition}[Generalized Maximum Principle]\label{Prop.GMP}
Let~$M$ be a linear elliptic second-order differential operator
with bounded coefficients in a bounded open set $U\subset\Real^d$;
the principal part of~$M$ is assumed to be formed
by a \emph{negative} definite matrix.
Suppose that there are two functions
$
  u, w \in \Smooth^2(U)\cap\Smooth^0(\overline{U})
$
satisfying:
$$
  \begin{aligned}
    M u &\leq 0 \quad &\mbox{in}& \quad U \,,
    \qquad &
    M w &\geq 0 \quad &\mbox{in}& \quad U \,,
    \\
    u &\geq 0 \quad &\mbox{in}& \quad U \,,
    \qquad &
    w &> 0 \quad &\mbox{in}& \quad \overline{U} \,.
  \end{aligned}
$$
Then
$$
  \sup_U(u/w) \, \leq \, \sup_{\partial U}(u/w) \,.
$$
\end{Proposition}
\noindent
To prove this, observe that the quotient $u/w$
is a subsolution of an elliptic operator
for which the usual maximum principle holds true
(\cf~\cite[Sec.~2.5]{PW}).

We shall also need two simple lemmata.
Writing $\Ball(t,r)$ for an open $(d-1)$-di\-men\-sional ball
of radius~$r>0$ centred at~$t\in\Real^{d-1}$
and abbreviating $\Ball_r:=\Ball(0,r)$,
we denote by~$\nu_1(r)$ the first eigenvalue
of the Dirichlet Laplacian in $\Ball_r$.
\begin{Lemma}\label{Lem.Ball}
Given a positive constant~$C$,
let $r>0$ be sufficiently small so that $\nu_1(r)>C$.
Then for any $f\in\Smooth^0(\partial\Ball_r)$,
the boundary problem
\begin{equation*}
  \left\{
  \begin{aligned}
  -\Delta v - C v &= 0
  \quad &\mbox{in}& \quad \Ball_r
  \,,
  \\
  v &= f
  \quad &\mbox{on}& \quad \partial\Ball_r
  \,,
  \end{aligned}
  \right.
\end{equation*}
has a unique solution
$
  v \in \Smooth^0(\overline{\Ball_r}) \cap \Smooth^2(\Ball_r)
$
and
$$
  v(0) = \alpha \ \frac{1}{|\partial\Ball_r|}
  \int_{\partial\Ball_r} f
$$
with some positive constant
$
  \alpha = \alpha(r,C,d)
$
such that
$
  r \mapsto \alpha(r,C,d)
$
is increasing.
Furthermore, $f>0$ implies $\inf v>0$.
\end{Lemma}
\begin{Lemma}\label{Lem.Fubini}
Let $u\in\sii(\Real^{d-1})$
be such that $\supp u \subseteq \overline{\cross}$,
and let $\delta>0$.
Then for any $t \in \overline{\cross}$,
there exists
$
  r = r(t,u,\delta) \in (0,\delta]
$
such that
$$
  \frac{1}{|\partial\Ball_r|}
  \int_{\partial\Ball(t,r)} |u|
  \ \leq \
  \frac{1}{\ |\Ball_{\delta}|^{1/2}}
  \, \|u\|_{\sii(\cross)}
  \,.
$$
\end{Lemma}
Lemma~\ref{Lem.Ball} can be established by standard arguments,
using the positive and rotationally symmetric eigenfunction
of the Dirichlet Laplacian in a larger ball
(\cf~also~\cite[Rem.~1]{Jerison}),
while Lemma~\ref{Lem.Fubini} follows easily by Fubini's theorem.

Now we are ready to prove
\begin{Theorem}\label{Thm.C0c}
Suppose Assumption~\ref{Ass.conv} holds.
For any $N\in\NatStar$, there exist
positive constants
$
  \eps_0 = \eps_0(N,L,C_{\ref{Ass.conv}},\cross)
$
and
$
  C = C(N,L,C_{\ref{Ass.conv}},\cross,d)
$
such that for all $\eps\leq\eps_0$,
$$
  \forall n\in\{1,\dots,N\}, \qquad
  \|\psi_n-\psi^0_n\|_{\Smooth^0(\overline{\cyl})}
  \ \leq \ C\,\eps
  \,.
$$
\end{Theorem}
\begin{proof}
Fix $n\in\NatStar$,
assume that~$\eps$ is sufficiently small so that the conclusion
of Theorem~\ref{Thm.H2c} holds true
with a positive constant~$C_1$
and recall the definition~(\ref{psi}).
Defining the operator
$$
  M :=
  -\partial_1 a_\eps \partial_1
  + \eps^{-2} \left(-\Delta'-E_1\right)
  + b_\eps
  \,,
$$
the first equation of~(\ref{identity})
together with~(\ref{Ass.conv'}) yields
$$
  M(\psi^\pm+\eps)
  = \pm f_\eps
  - (\eps^{-1} E_1 - b_\eps\eps)
  \leq 0
  \qquad\mbox{in}\quad
  \cyl^\pm := \{(s,t)\in\cyl\,|\ \psi^\pm > 0 \}
$$
provided $\eps^2 \leq E_1/(2C_A')$,
where
$
  \psi^\pm := \max\{\pm \psi,0\}
$.
That is, $\psi^\pm+\eps$ is a positive subsolution of~$M$ in~$\cyl^\pm$.
Our strategy will be to find a supersolution
appropriate for the comparison argument
included in Proposition~\ref{Prop.GMP}.
Defining
$$
  u(t) :=
  \begin{cases}
    \sup_{s \in I} |\psi(s,t)|
    & \quad\mbox{if}\quad t\in\overline{\cross}
    \,,
    \\
    0
    & \quad\mbox{if}\quad t\in\Real^{d-1}\setminus\overline{\cross}
    \,,
  \end{cases}
$$
we obviously have
\begin{equation}\label{u-H1}
  \|u\|_{\sii(\cross)}
  \leq L^{1/2} \, \|\partial_1\psi\|_{\Hilbert}
  \leq L^{1/2} \, C_1 \, \eps
  \,.
\end{equation}
Let $\delta>0$ be (uniquely) determined by the condition
$\nu_1(\delta)=3 E_1$.
Fix
$
  t_0 \in \overline{\cross}
$,
and let $r=r(t_0,u,\delta)$ be the corresponding radius
determined by Lemma~\ref{Lem.Fubini}.
Finally, define a function $w$ on $I\times\Ball(t_0,r)$
by putting $w:=\id \otimes v$,
where~$v$ is the solution to
\begin{equation*}
  \left\{
  \begin{aligned}
  -\Delta v - 2 E_1 v &= 0
  \quad &\mbox{in}& \quad \Ball(t_0,r)
  \,,
  \\
  v &= u + \eps
  \quad &\mbox{on}& \quad \partial\Ball(t_0,r)
  \,.
  \end{aligned}
  \right.
\end{equation*}
By~Lemma~\ref{Lem.Ball}, $v$ indeed exists
and it is continuous and positive up to the boundary.
Moreover, the minimum of~$v$ is achieved
on the boundary because~$v$ is superharmonic.
Since
$$
  M w = (\eps^{-2} E_1 + b_\eps) w
  \geq 0
  \qquad\mbox{in}\quad
  I \times \Ball(t_0,r)
$$
provided
$
  \eps^2 \leq E_1/C_{\ref{Ass.conv}}'
$,
and since~$\psi^\pm$ is either equal to zero
or not greater than~$w-\eps$ on the boundary of
$
  U^\pm := \cyl^\pm \cap [I \times \Ball(t_0,r)]
$,
Proposition~\ref{Prop.GMP} yields
$$
  \forall (s,t) \in \overline{U^\pm}, \qquad
  \psi^\pm(s,t) + \eps
  \leq v(t)
  \,.
$$
In particular, Lemma~\ref{Lem.Ball} gives
\begin{equation*}
  \forall s \in \overline{I}, \qquad
  |\psi(s,t_0)| + \eps
  \leq v(t_0)
  = \alpha(r,2E_1,d)
  \ \frac{1}{|\partial\Ball_r|}
  \int_{\partial\Ball(t_0,r)} (u+\eps)
  \,,
\end{equation*}
where~$\alpha$ is the coefficient of Lemma~\ref{Lem.Ball}.
Using Lemmata~\ref{Lem.Ball}, \ref{Lem.Fubini} and~(\ref{u-H1}),
we conclude that
$$
  \forall (s,t_0) \in \overline{\cyl}, \qquad
  |\psi(s,t_0)|
  \leq
  \alpha(\delta,2E_1,d)
  \, \big(
  1 + |\Ball_{\delta}|^{-1/2} \, L^{1/2} \, C_1
  \big)
  \, \eps
$$
because~$t_0$ was chosen arbitrarily.
\end{proof}
%

\subsection{\mbox{$\Smooth^{0}$-convergence of transverse
    derivatives on the boundary}}\label{Sec.Lip}
%
Now we use the maximum principle
to derive a Lipschitz-type condition for~$\psi$,
which will play a crucial role in our proof
of convergence of nodal sets.
\begin{Theorem}\label{Thm.Lip}
Suppose Assumption~\ref{Ass.conv} holds.
For any $N\in\NatStar$, there exist
positive constants
$
  \eps_0 = \eps_0(N,L,C_{\ref{Ass.conv}},\cross)
$
and
$
  C = C(N,L,C_{\ref{Ass.conv}},\cross,d)
$
such that for all $\eps\leq\eps_0$,
$$
  \forall n\in\{1,\dots,N\}, \qquad
  \forall t_0\in \partial\cross, \ \,
  \sup_{(s,t) \in I\times\cross}
  \left|
  \frac{\psi_n(s,t)}{|t-t_0|}
  - \frac{\psi^0_n(s,t)}{|t-t_0|}
  \right|
  \ \leq \ C\,\eps
  \,.
$$
\end{Theorem}
\begin{proof}
Rewrite the first line of~(\ref{identity}) as
\begin{equation*}
  M\psi :=
  -\eps^2\partial_1\!\left(a_\eps\,\partial_1 \psi\right)
  -\Delta'\psi
  = \eps^2 f_\eps - \eps^2 b_\eps\psi + E_1 \psi
  =: F_\eps
  \,.
\end{equation*}
Theorem~\ref{Thm.C0c} and~(\ref{Ass.conv'}) imply that
for all $\eps\leq\eps_1$,
$$
  |M\psi|
  \leq \|F_\eps\|_{\Smooth^0(\overline{\cyl})}
  \leq C_1 \eps
  \qquad\mbox{in}\quad \cyl
$$
with some positive constants
$
  \eps_1 = \eps_1(n,L,C_{\ref{Ass.conv}},\cross)
$
and
$
  C_1 = C_1(n,L,C_{\ref{Ass.conv}},\cross,d)
$.
We are now inspired by \cite[Prob.~3.6]{Gilbarg-Trudinger}
to construct a supersolution suitable for a comparison argument.
Since the boundary~$\partial\cross$ is of class $\Smooth^\infty$,
there exists a positive number $r=r(\cross)$ such that
for every boundary point $t_0\in\partial\cross$,
there is an exterior point $\tau\in\Real^{d-1}\setminus\overline{\cross}$
satisfying
$
  \overline{\Ball(\tau,r)} \cap \overline{\cross} = \{t_0\}
$.
Recall also the definition~(\ref{def.a}) of $a$.
In the cylindrical layer
$$
  U :=
  I \times
  \big[
  \Ball(\tau,r+2 a)
  \setminus\overline{\Ball(\tau,r)}
  \big]
  \supset \cyl
  \,,
$$
consider a (positive) function~$w$ defined by
$$
  w := 1 \otimes v \,,
  \qquad
  v(t) := \beta \,
  \big( e^{-\alpha r} - e^{-\alpha |t-\tau|} \big)
  \,,
$$
where~$\alpha$ and~$\beta$ are positive parameters yet to be determined.
Direct calculation gives for $(s,t) \in U$,
\begin{equation*}
  (M w)(s,t)
  = \beta\,\alpha\, e^{-\alpha |t-\tau|}
  \big[
  \alpha - (d-2) \, |t-\tau|^{-1}
  \big]
  \geq
  \beta\,\alpha\, e^{-\alpha (r+2a)}
  \big[
  \alpha - (d-2)\,r^{-1}
  \big]
  \,,
\end{equation*}
where the inequality hold
provided $\alpha \geq (d-2) \, r^{-1}$.
Choosing, \eg,
$$
  \alpha := (d-1) \, r^{-1}
  \qquad\mbox{and}\qquad
  \beta :=
  \alpha^{-1} \, r \, e^{\alpha\,(r+2a)} \,
  \|F_\eps\|_{\Smooth^0(\overline{\cyl})}
  \,,
$$
we therefore have
$$
  M w
  \geq \|F_\eps\|_{\Smooth^0(\overline{\cyl})}
  \qquad\mbox{in}\quad U
  \,.
$$
Since $\psi=0$ on~$\partial\cyl$
and~$w$ is non-negative on~$\partial U$,
the maximum principle (\cf~\cite[Thm.~3.3]{Gilbarg-Trudinger}) yields
$
  |\psi| \leq w
$
in the closure of $\cyl$.
However, for $(s,t) \in U$,
$$
  w(s,t) \leq \alpha \, \beta \, e^{-\alpha r} \,
  \big(|t-\tau|-|\tau-t_0|\big)
  \leq \alpha \, \beta \, e^{-\alpha r} \,
  |t-t_0|
  \,,
$$
and the claim of Theorem therefore follows with
$
  C := C_1 \, r e^{(d-1)r^{-1} 2 a}
$.
\end{proof}

The part~(ii) of Theorem~\ref{Thm.main.intro}
follows as a consequence of Theorem~\ref{Thm.Lip}.

\subsection{Convergence of nodal domains}\label{Sec.Nodal}
%
Now we are in a position to establish
the main result of this section,
referring to Section~\ref{Sec.longitudinal}
for the definition of the one-dimensional eigenfunction~$\phi_n$
and the discussion of its nodal set.
\begin{Theorem}\label{Thm.main}
Suppose Assumption~\ref{Ass.conv} holds.
For any $N\in\NatStar$,
there exist positive constants
$
  \eps_0 = \eps_0(N,L,C_{\ref{Ass.conv}},\cross)
$
and
$
  C = C(N,L,C_{\ref{Ass.conv}},\cross,d)
$
such that for all $\eps\leq\eps_0$
and any $n\in\{1,\dots,N\}$,
\begin{enumerate}
\item[(i)]
$
  \forall (s,t)\in\cyl, \quad
  (s,t) \in \nodal(\psi_n)
  \ \Longrightarrow \
  \dist\big(s,\nodal(\phi_n)\big)
  \, \leq \, C\,\eps
$\,,
\item[(ii)]
$
  \forall (s,t)\in\cyl, \quad
  \dist\big(s,\nodal(\phi_n)\big) > C\,\eps
  \ \Longrightarrow \
  \sgn\psi_n(s,t)
  = \sgn\phi_n(s)
$\,.
\end{enumerate}
\end{Theorem}
\begin{proof}
Fix $n\in\Nat$, $n \geq 2$,
and recall the definitions of~$I_n$, $I_n^k$ and $s_0^k(n)$
introduced in Section~\ref{Sec.longitudinal}.
Let~$\eps$ be sufficiently small so that conclusions
of Proposition~\ref{Prop.nodal}.(ii) and Theorem~\ref{Thm.Lip}
hold true with positive constants~$c_1$ and~$C_1$, respectively.
Combining Theorem~\ref{Thm.Lip}, Proposition~\ref{Prop.nodal}.(ii)
and~(\ref{ef}), we have for all $(s,t) \in I_n\times\cross$,
\begin{equation}\label{positivity}
  \frac{\psi_n(s,t)}{\dist(t,\partial\cross)}
  \, \sgn\phi_n(s)
  \ \geq \ c_1 \, \dist(s,\partial I_n)
  - C_1\,\eps
  \ > \ 0
\end{equation}
provided
$
  \dist(s,\partial I_n)
  > C  \, \eps
$,
with
$
  C  := c_1^{-1}  C_1
$.
This establishes~(i) and~(ii) with $\nodal(\phi_n)$
being replaced by
$
  \partial I_n = \nodal(\phi_n) \cup \{0,L\}
$.

It remains to show that
$
  \nodal(\psi_n) \cap (I_n'\times\cross) = \varnothing
$,
where
$$
  I_n' :=
  \big(0,s_0^1(n)-C  \, \eps\big)
  \cup
  \big(s_0^{n-1}(n)+C  \, \eps,L\big)
  \,.
$$
Let~$c_2$ be the positive constant
determined by Proposition~\ref{Prop.Sturm}.
If
$
  \eps \leq c_2/(4C)
$
so that by Proposition~\ref{Prop.Sturm}.(iv)
there does exist an open non-empty interval
$J_n^k \subset I_n^k$ satisfying
$
  \dist(J_n^k,\partial I_n)
  > C \, \eps
$
for all $k\in\{0,\dots,n-1\}$,
(\ref{positivity}) also implies that each $J_n^k\times\cross$
belongs to one distinct nodal domain of~$\psi_n$
(\ie\ connected component of $\cyl\setminus\nodal(\psi_n)$).
In particular, $\psi_n$ has already~$n$ nodal domains,
and Courant's nodal domain theorem~\cite{Alessandrini_1998}
implies that~$\psi_n$ cannot change sign in each of
the connected components of $I_n'\times\cross$.
But then~$\psi_n$ cannot vanish in $I_n'\times\cross$,
\eg, because of the Harnack inequality
(\cf~\cite[Thm.~8.20]{Gilbarg-Trudinger}).
\end{proof}

This concludes the proof of Theorem~\ref{Thm.main.intro}.
Corollary~\ref{Corol.main} is a consequence
of the property~(iii) of Theorem~\ref{Thm.main.intro}
and Courant's nodal domain theorem~\cite{Alessandrini_1998}.

\setcounter{equation}{0}
\section{Curved tubes}\label{Sec.tubes}
%
In this section, we consider the Dirichlet Laplacian
in curved tubes of shrinking cross-section.
Using special curvilinear coordinates,
we transform the Laplacian in a tube to a unitarily equivalent
Schr\"odinger-type operator of the form~(\ref{operator})
in a straight tube, and apply Theorem~\ref{Thm.main.intro}.
The  necessary geometric preliminaries follow
the lines of Section~\ref{Sec.tubes.intro},
but we proceed in more details.

\subsection{The reference curve}
%
Let us precise what we meant by the \emph{appropriate}
Frenet frame of the reference curve~$\curve$
in the beginning of Section~\ref{Sec.tubes.intro}:
\begin{Assumption}\label{Ass.Frenet}
$\curve$ possesses a positively oriented
Frenet frame $\{e_1,\dots,e_d\}$
with the properties that
\begin{enumerate}
\item[(i)]
$
  e_1 = \dot{\curve}
$\,,
\item[(ii)]
$
  \forall i \in \{1,\dots,d\},\quad
  e_{i} \in \Smooth^\infty(\overline{I};\Real^d)
$\,,
\item[(iii)]
$
  \forall i \in \{1,\dots,d-1\}, \
  \forall s\in\overline{I},
  \quad
  \dot{e}_i(s)
$
lies in the span of $e_1(s),\dots,e_{i+1}(s)$\,.
\end{enumerate}
\end{Assumption}
\noindent
Here and in the sequel, the \emph{dot} denotes the derivative.
\begin{Remark}\label{Rem.Frenet}
Recall~\cite[Sec.~1.2]{Kli} that a Frenet frame
is by definition a moving (orthonormal) frame
such that for all $i\in\{1,\dots,d\}$ and $s \in \overline{I}$,
$\Gamma^{(i)}(s)$ lies in the span of $e_1(s),\dots,e_i(s)$.
A sufficient condition to ensure the existence
of the Frenet frame of Assumption~\ref{Ass.Frenet}
is to require that for all $s \in \overline{I}$,
the vectors
$
  \dot{\curve}(s), \curve^{(2)}(s), \dots, \curve^{(d-1)}(s)
$
are linearly independent (\cf~\cite[Prop.~1.2.2]{Kli}).
This is always satisfied if $d=2$.
However, we do not assume \emph{a priori}
this non-degeneracy condition for $d \geq 3$
because it excludes, \eg, the curves such that
$\curve(J)$ is a straight segment
for some subinterval $J \subseteq I$.
\end{Remark}

Let
$
  \curv\equiv(\curv_{ij})_{i,j=1}^d
$
be the matrix-valued function
defined by the (Serret-Frenet) formulae
\begin{equation}\label{Frenet}
  \dot{e}_i = \sum_{j=1}^d \curv_{ij} \, e_j
  \,, \qquad
  i \in \{1,\dots,d\}
  \,.
\end{equation}
By virtue of Assumption~\ref{Ass.Frenet},
$\curv$ has a skew symmetry
\begin{equation*}
  \curv_{ij}=-\curv_{ji}
  \,, \qquad
  i,j \in \{1,\dots,d\}
\end{equation*}
and $\curv_{ij}=0$ for $j>i+1$.
We define by $\kappa_i := \curv_{i,i+1}$
the $i^\mathrm{th}$ curvature of~$\curve$.
It will be also convenient to introduce the submatrix
$
  \curv' := (\curv_{\mu\nu})_{\mu,\nu=2}^d
$.

\subsection{The Tang frame}\label{Sec.Tang}
%
We introduce now another moving frame along~$\curve$,
which better reflects the geometry of the curve.
We shall refer to it as the \emph{Tang frame}
because it is a natural generalization
of the Tang frame known from the theory
of three-dimensional waveguides~\cite{Tsao-Gambling_1989}.
Our construction follows the generalization
introduced in~\cite{ChDFK}.

Let
$
  \rot' \equiv
  (\rot_{\mu\nu})_{\mu,\nu=2}^d
$
be the $(d-1)\times(d-1)$ matrix-valued function
defined by the initial-value problem
\begin{equation}\label{diff.eq}
  \left\{
  \begin{aligned}
    \dot{\rot}' + \rot' \, \curv' &= 0
    &\quad\mbox{in}& \quad [0,L]
    \,,
    \\
    \rot' &= \rot_0'
    &\quad\mbox{at}& \quad 0
    \,,
  \end{aligned}
  \right.
\end{equation}
where $\rot_0'$ is a rotation matrix in $\Real^{d-1}$, \ie,
\begin{equation}\label{rotation}
  \det(\rot_0') = 1
  \,, \qquad
  \rot_0' \, {\rot_0'}^{\!\mathrm{T}} = \id \,.
\end{equation}
Here~``$\mathrm{T}$" denotes the transpose operation and~$1$ stands
both for a scalar number and an identity matrix.
The solution of~(\ref{diff.eq}) exists,
it is unique (for a given~$\rot_0'$),
belongs to
$
  \Smooth^\infty(\overline{I};\Real^{(d-1)^2})
$
and satisfies the conditions~(\ref{rotation})
in all~$\overline{I}$
(\cf~\cite[Sec.~2.2]{ChDFK} for more details).

We extend~$\rot'$ to a $d\times d$ matrix-valued function
$
  \rot \equiv
  (\rot_{ij})_{i,j=1}^d
$
by setting
\begin{equation}
  \rot
  :=
  \begin{pmatrix}
    1 & 0 \\
    0 & \rot'
  \end{pmatrix} ,
\end{equation}
and introduce the Tang frame
$
  \{\tilde{e}_1, \dots, \tilde{e}_d\}
$
by
\begin{equation}\label{Tang}
  \tilde{e}_i := \sum_{j=1}^d \rot_{ij} \, e_j
  \,, \qquad
  i\in\{1,\dots,d\}
  \,.
\end{equation}
Combining~(\ref{diff.eq}) and~(\ref{Frenet})
together with the properties of~$\curv$,
one easily finds
\begin{equation}\label{Frenet.bis}
  \dot{\tilde{e}}_1 = \kappa_1 \, e_2
  \qquad\textrm{and}\qquad
  \dot{\tilde{e}}_\mu = - \kappa_1 \, \rot_{\mu 2} \, e_1
  \,, \qquad
  \mu\in\{2,\dots,d\}
  \,.
\end{equation}
%

\subsection{Tubes}
%

As in Section~\ref{Sec.tubes.intro},
we introduce a mapping~$\tubemap$
from a straight tube~$\cyl$ to~$\Real^d$
by~(\ref{tube.map.intro}).
Recalling the relation~(\ref{Tang})
between the Frenet and Tang frames,
it is clear that the image $\tubemap(\cyl)$
is obtained by ``translating'' the cross-section~$\eps\cross$
along the curve~$\curve$ in a special way
(it ``rotates'' with respect to the Tang frame).

Obviously,
$
  \tubemap\in\Smooth^\infty(\overline{\cyl};\Real^d)
$
for all $\eps>0$.
Furthermore, $\tubemap$ is an immersion
provided~$\eps$ is small enough.
This can be seen as follows.
Let
$
  G\equiv(G_{ij})_{i,j=1}^d
$
be the metric induced by the mapping~$\tubemap$, \ie,
$
  G_{ij} :=
  (\partial_i\tubemap) \cdot (\partial_j\tubemap)
$,
where ``$\cdot$" denotes the inner product in~$\Real^d$.
Using the orthonormality of the Tang frame,
relations~(\ref{Frenet.bis})
and Assumption~\ref{Ass.Frenet}.(i),
one easily establishes the formulae
~(\ref{metric}) and~(\ref{Jacobian}).
Consequently,
$$
  |G| := \det(G) = \eps^{2(d-1)} \, h^2
  \,.
$$
Since~$\rot$ satisfies the orthogonality condition~(\ref{rotation})
and $|t| \leq a$
(recall that~$\cross$ is assumed to have
its centre of mass at the origin of~$\Real^{d-1}$
and~$a$ is given by~(\ref{def.a})),
the restriction~(\ref{class}) yields
\begin{equation}\label{ellipticity}
  1 - C_\curve \, a \, \eps
  \ \leq \ h \ \leq \
  1 + C_\curve \, a \, \eps
  \,.
\end{equation}
In particular, $h$ does not vanish in~$\cyl$ provided
\begin{equation}\label{Ass.basic}
  \eps < (C_\curve\,a)^{-1}
  \,,
\end{equation}
and it follows by the inverse function theorem
that~$\tubemap$ induces a local $\Smooth^\infty$-diffeomorphism.
This shows that~$\tubemap$ is an immersion
for all positive~$\eps$ satisfying~(\ref{Ass.basic}).

In Definition~\ref{Def.tubes},
we have introduced the tube~$\tube$
to be the manifold~$\cyl$ equipped
with the Riemannian metric~$G$.
The symbol~$d\vol$ will denote the volume measure on~$\tube$,
\ie,
$$
  d\vol := |G(s,t)|^{1/2} \, ds \, dt
  = \eps^{d-1} \, h(s,t)\,ds\,dt
  \,,
$$
where~$dt \equiv dt_2 \dots dt_d$
stays for the $(d-1)$-dimensional Lebesgue measure on~$\cross$.

\subsection{The Laplacian}\label{Sec.Laplacian}
%
It is a general fact that the Laplacian~$-\Delta_G$
(as a differential expression)
in a manifold~$\cyl$ equipped with the metric~$G$
can be written as
\begin{equation*}
  -\Delta_G
  \ = \
  -|G|^{-1/2} \sum_{i,j=1}^d
  \partial_i |G|^{1/2} G^{ij} \partial_j
  \,,
\end{equation*}
where~$G^{ij}$ denote the coefficients
of the inverse matrix~$G^{-1}$.
We introduce the \emph{Dirichlet Laplacian}~$-\Delta_D^\tube$
(as a differential operator) in the tube $\tube\equiv(\cyl,G)$
to be the operator in the Hilbert space
$
  \sii(\tube) \equiv \sii(\cyl,d\vol)
$
defined by
\begin{equation}\label{Laplace}
  -\Delta_D^\tube\psi
  := -\Delta_G\psi
  \,, \qquad
  \psi \in \Dom(-\Delta_D^\tube)
  := \Hilbert_0^1(\tube) \cap \Hilbert^2(\tube)
  \,.
\end{equation}
It can be verified directly that~$-\Delta_D^\tube$ is self-adjoint.
Alternatively, this can be deduced
from the following unitary equivalence
we shall need anyway.

Let~$\mathcal{U}$ be the unitary operator defined by
\begin{equation}\label{unitary}
  \mathcal{U}: \sii(\tube) \to \sii(\cyl) :
  \big\{ \psi \mapsto  |G|^{1/4}\psi \big\}
  \,.
\end{equation}
Setting $H:=\mathcal{U}\,(-\Delta_D^\tube)\,\mathcal{U}^{-1}$,
one can check that
\begin{equation}\label{Hamiltonian.potential}
  H = -\partial_1 h^{-2} \partial_1
  - \eps^{-2} \Delta' + V
  \,, \qquad
  \Dom(H) = \Hilbert_2
  \,,
\end{equation}
where the space~$\Hilbert_2$ is defined in~(\ref{domain}) and
\begin{equation}\label{potential}
  V \ := \
  - \frac{5}{4} \, \frac{(\partial_1 h)^2}{h^4}
  + \demi \, \frac{\partial_1^2h}{h^3}
  - \frac{1}{4} \, \frac{|\nabla'h|^2}{\eps^2 h^2}
  + \demi \, \frac{\Delta' h}{\eps^2 h}
  \,.
\end{equation}

Actually, (\ref{Hamiltonian.potential}) with~(\ref{potential})
is a general formula valid for any $\Smooth^\infty$-smooth metric
of the form~(\ref{metric}).
In our special case when~$h$ is given by~(\ref{Jacobian}),
we find easily that
$\partial_\mu h=-\eps\,\kappa_1 \, \rot_{\mu 2}$
and $\partial_\mu\partial_\nu h=0$
for any $\mu,\nu\in\{2,\dots,d\}$,
and therefore
\begin{equation}\label{potential.bis}
  V \ = \
  - \frac{1}{4} \, \frac{\kappa_1^2}{h^2}
  + \demi \, \frac{\partial_1^2 h}{h^3}
  - \frac{5}{4} \, \frac{(\partial_1 h)^2}{h^4}
  \,.
\end{equation}
Moreover, (\ref{diff.eq})~gives
\begin{equation*}
\begin{aligned}
  \partial_1 h(\cdot,t)
  &\ = \ \eps \sum_{\mu,\nu,\rho=2}^d
  t_\mu\,\rot_{\mu\nu}
  \big(\dot{\curv}_{\nu 1}
  -\curv_{\nu\rho}\curv_{\rho 1}\big)\,,
  \\
  \partial_1^2 h(\cdot,t)
  &\ = \ \eps \sum_{\mu,\nu,\rho,\sigma=2}^d
  t_\mu\,\rot_{\mu\nu}
  \big(\ddot{\curv}_{\nu 1}
  -\dot{\curv}_{\nu\rho}\curv_{\rho 1}
  -2\;\!\curv_{\nu\rho}\dot{\curv}_{\rho 1}
  +\curv_{\nu\rho}\curv_{\rho\sigma}\curv_{\sigma 1}
  \big)
  \,.
\end{aligned}
\end{equation*}

Using~(\ref{ellipticity}), it is easy to see that~$H$
is uniformly elliptic
with uniformly $\Smooth^\infty$-smooth coefficients
for all~$\eps$ satisfying~(\ref{Ass.basic}).
Consequently, $H$ (and therefore~$-\Delta_D^\tube$)
is self-adjoint by the same reasoning
as in Section~\ref{Sec.operator}.

\subsection{Thin tubes}\label{Sec.application}
%
It remains to apply the results of Section~\ref{Sec.Conv},
namely Theorem~\ref{Thm.main.intro}, to~$H$.
Let us define the operator
\begin{equation}\label{operator.bis}
  T := H - \eps^{-2} E_1
  \,,
\end{equation}
which is indeed of the form~(\ref{operator}),
with $a_\eps:=h^{-2}$ and $V_\eps:=V$.
It is important to notice that,
while the eigenvalues of~$T$ are just eigenvalues of~$H$
shifted by~$-\eps^{-2} E_1$,
the operators have in fact the same eigenfunctions.

Let us now verify Assumption~\ref{Ass.conv} for~(\ref{operator.bis}).
Assuming~(\ref{Ass.basic}),
the functions $h^{-2}$ and $V$
obviously belong to $\Smooth^\infty(\overline{\cyl})$
(they are even analytic in~$t$),
since it is true for~$\tubemap$
(as a consequence of the geometric Assumption~\ref{Ass.Frenet});
this checks Assumption~\ref{Ass.conv}.(i).
In view of~(\ref{ellipticity}) and~(\ref{Ass.basic}),
Assumption~\ref{Ass.conv}.(ii) holds true
with $C_{\ref{Ass.conv}}=4$.
Strengthen now~(\ref{Ass.basic}) to
$
  \eps\leq(2C_\curve a)^{-1}
$,
so that $h \geq 1/2$ uniformly in~$\eps$
by virtue of~(\ref{ellipticity}).
Since
\begin{align*}
  \|1-h^2\|_{\Smooth^0(\overline{\cyl})}
  &\leq \eps \, a \, \|\kappa_1\|_{\Smooth^0(\overline{I})}
  \big(
  2+\eps\, a \, \|\kappa_1\|_{\Smooth^0(\overline{I})}
  \big)
  \,,
  \\
  \|\partial_1 h\|_{\Smooth^0(\overline{\cyl})}
  &\leq \eps \, a \,
  \big(
  \|\dot{\kappa}_1\|_{\Smooth^0(\overline{I})}
  + \|\curv\|_{\Smooth^0(\overline{I};\Real^{d^2})}^2
  \big)
  \,,
  \\
  \|\partial_1^2 h\|_{\Smooth^0(\overline{\cyl})}
  &\leq \eps \, a\,
  \big(
  \|\ddot{\kappa}_1\|_{\Smooth^0(\overline{I})}
  + 3 \|\dot{\curv}\|_{\Smooth^0(\overline{I};\Real^{d^2})}
  \|\curv\|_{\Smooth^0(\overline{I};\Real^{d^2})}
  + \|\curv\|_{\Smooth^0(\overline{I};\Real^{d^2})}^3
  \big)
  \,,
  \\
  \|\partial_\mu h\|_{\Smooth^0(\overline{\cyl})}
  &\leq \eps \, a \, \|\kappa_1\|_{\Smooth^0(\overline{I})}
  \,,\qquad \mu\in\{2,\dots,d\} \,,
\end{align*}
it is easy to see from~(\ref{potential}) by applying~(\ref{class})
that Assumptions~\ref{Ass.conv}.(iii) and~(iv)
hold true with some positive~$C_{\ref{Ass.conv}}$
depending on~$C_\curve$, $a$ and~$d$,
and with~$v_0$ as defined in~(\ref{potential.tube}).

The application of Theorem~\ref{Thm.main.intro}
to the operator~(\ref{operator.bis}) is now straightforward.
In order to deduce from it Theorem~\ref{Thm.tube.intro},
it remains to recall that~$H$ and~$-\Delta_D^\tube$
have the same eigenvalues
and their eigenfunctions are related by (\ref{unitary}).

\setcounter{equation}{0}
\section{Strips on surfaces}\label{Sec.Strips}
%
As another application of Theorem~\ref{Thm.main.intro},
let us consider the situation where the ambient space
of the tube is a general Riemannian manifold
instead of the Euclidean space~$\Real^d$.
We restrict ourselves to the case $d=2$, \ie,
the tube~$\tube$ will be a strip about
a curve in an (abstract) two-dimensional surface.
We refer to~\cite{K1,KT,K2} for geometric details
and basic spectral properties of~$-\Delta_D^\tube$
in the infinite case $I=\Real$.

Consider a $\Smooth^\infty$-smooth
connected complete non-compact two-dimensional
Riemannian manifold~$\ambient$ of bounded Gauss curvature~$K$,
and a $\Smooth^\infty$-smooth curve $\curve:\overline{I}\to\ambient$
which is assumed to be parametrized by arc length and embedded.
Let $N$ be the unit normal vector field along~$\curve$,
which is uniquely determined as the $\Smooth^\infty$-smooth mapping
from~$\overline{I}$ to the tangent bundle of~$\ambient$
by requiring that
$N(s)$ is orthogonal to the derivative~$\dot{\curve}(s)$
and that $\{\dot{\curve}(s),N(s)\}$ is positively oriented
for all $s\in\overline{I}$ (\cf~\cite[Sec.~7.B]{Spivak4}).
We denote by~$\kappa$ the corresponding curvature of~$\curve$
defined by the Frenet formula
and note that its sign is uniquely determined
up to the re-parametrization of~$\curve$
($\kappa$ is the geodesic curvature of~$\curve$
if~$\ambient$ is embedded in~$\Real^3$).

Without loss of generality,
$\cross$~can be chosen as $(-1,1)$.
For sufficiently small positive~$\eps$,
we define a mapping~$\tubemap$
from~$\cyl$ to~$\ambient$ by setting
\begin{equation}\label{exponential}
  \tubemap(s,t) := \exp_{\curve(s)}\big(\eps\,t\,N(s)\big)
  \,,
\end{equation}
where~$\exp_p$ is the exponential map of~$\ambient$
at $p\in\ambient$.
Note that $s\mapsto\tubemap(s,t)$ traces the curves
parallel to~$\curve$ at a fixed distance~$\eps |t|$,
while the curve~$t\mapsto\tubemap(s,t)$
is a geodesic orthogonal to~$\curve$
for any fixed~$s$.
Since~$\curve$ is compact, $\tubemap$ induces a diffeomorphism
of~$\cyl$ onto the image~$\tubemap(\cyl)$
provided~$\eps$ is small enough
(\cf~\cite[Sec.~3.1]{Gray}).
Consequently, $\tubemap$~induces a Riemannian metric~$G$ on~$\cyl$,
and we define the \emph{strip}~$\tube$ about~$\curve$
to be the manifold $(\cyl,G)$.
It follows by the generalized Gauss lemma~\cite[Sec.~2.4]{Gray}
that the metric acquires the diagonal form:
\begin{equation*}
  G = \diag\big(h^2,\eps^2\big) \,,
\end{equation*}
where~$h$ is a uniformly $\Smooth^\infty$-smooth function on~$\cyl$
satisfying the Jacobi equation
\begin{equation}\label{Jacobi}
  \partial_2^2 h + \eps^2\,K\,h = 0
  \qquad\textrm{with}\qquad\left\{
  \begin{aligned}
    h(\cdot,0) &= 1 \,, \\
    \partial_2 h(\cdot,0) &= -\eps\,\kappa \,,
  \end{aligned}
  \right.
\end{equation}
where~$K$ is considered as a function
of the (Fermi) ``coordinates" $(s,t)$
determined by~(\ref{exponential}).

Since the metric~$G$ is of the form~(\ref{metric}),
the Dirichlet Laplacian in the manifold~$\tube$
(defined in the same manner as~(\ref{Laplace})),
is unitarily equivalent to the operator~$H$
defined by~(\ref{Hamiltonian.potential}) with~(\ref{potential}).
It is easy see that Theorem~\ref{Thm.main.intro}
applies to the shifted operator $H-\eps^{-2}E_1$.
In particular, using~(\ref{Jacobi}),
we verify that Assumption~\ref{Ass.conv} holds true
with some positive constant depending
on the norms $\|\kappa\|_{\Smooth^2(\overline{I})}$
and $\|K\|_{\Smooth^2(\overline{I})}$,
and with the function~$v_0$ given this time by
(\cf~\cite[Sec.~3]{K1}):
$$
  v_0(s) := - \frac{\kappa(s)^2}{4} - \frac{K(s,0)}{2}
  \,.
$$
The latter determines spectral properties
of the one-dimensional operator~$S$ introduced in~(\ref{op.geom}),
namely, the set of eigenvalues $\{\mu_n\}_{n=1}^\infty$
and corresponding eigenfunctions $\{\phi_n\}_{n=1}^\infty$.
\begin{Theorem}\label{Thm.strips}
Let~$\tube$ be the strip of width~$2\eps$ defined above
as a tubular neighbourhood about a curve~$\curve$ embedded
in a two-dimensional Riemannian manifold~$\ambient$.
For any integer $N \geq 1$,
there exist positive constants~$\eps_0$ and~$C$
depending on $N,L$ and the geometries of~$\curve$ and~$\ambient$
such that for all $\eps \leq \eps_0$,
the claims~(i)--(iii) of Theorem~\ref{Thm.tube.intro} hold true.
\end{Theorem}

\paragraph{Acknowledgements.}
The authors would like to thank Diogo Gomes
and Denis I.~Borisov for useful remarks.
Most of the results presented here were obtained
during the stay of one of the authors (D.K.)
at the Departamento de Matem\'atica,
Instituto Superior T\'ecnico, Lisboa, Portugal;
the author expresses his gratitude to the host.
The work was partially supported by FCT/\-POCTI/\-FEDER, Portugal.
The second author (D.K.) was also supported by
the Czech Academy of Sciences and its Grant Agency
within the projects IRP AV0Z10480505 and A100480501,
and by the project LC06002 of the Ministry of Education,
Youth and Sports of the Czech Republic.

%
\addcontentsline{toc}{section}{References}

\providecommand{\bysame}{\leavevmode\hbox to3em{\hrulefill}\thinspace}

\end{document}